\documentclass[12pt,a4paper,reqno]{amsart}
\usepackage{amsmath,amssymb}
\usepackage[dvipsnames]{xcolor}
\usepackage[colorlinks=true,urlcolor=blue,linkcolor=blue,citecolor=blue]{hyperref}

\numberwithin{equation}{section}
\theoremstyle{plain}
\newtheorem{thm}{Theorem}[section]
\newtheorem{lem}[thm]{Lemma}

\newtheorem{claim}[thm]{Claim}
\newtheorem{cor}[thm]{Corollary}

\theoremstyle{definition}

\newtheorem{exa}[thm]{Example}

\newcommand*{\proofclaimname}{Proof}
\newenvironment{proofclaim}[1][\proofclaimname]{\begin{proof}[#1]}{\end{proof}}

\newcommand{\A}{\mathbf{A}}
\newcommand{\AB}{\mathbf{A^B}}
\newcommand{\ABxe}{(\mathbf{A^B})^{x_1,\dots,x_n}_{e_1,\dots,e_n}}

\newcommand{\B}{\mathbf{B}}
\newcommand{\C}{\mathbf{C}}
\newcommand{\D}{\mathbf{D}}
\newcommand{\F}{\mathbf{F}}

\renewcommand{\H}{\mathbf{H}}
\newcommand{\id}{\mathrm{id}}
\newcommand{\K}{\mathfrak{K}}
\newcommand{\R}{\mathbf{R}}
\newcommand{\N}{\mathbb{N}}

\newcommand{\Z}{\mathbb{Z}}

\DeclareMathOperator{\Aut}{Aut}
\DeclareMathOperator{\Homeo}{Homeo}

\DeclareMathOperator{\supp}{supp}

\DeclareMathOperator{\cl}{cl}

\newcommand{\st}{\: :\:}

\usepackage{changes}
\makeatletter
\AddToHook{cmd/added/before}{\def\Changes@AuthorColor{blue}}
\AddToHook{cmd/deleted/before}{\def\Changes@AuthorColor{red}}
\AddToHook{cmd/replaced/before}{\def\Changes@AuthorColor{blue}}
\makeatother

\makeatletter
\@namedef{subjclassname@2020}{\textup{2020} Mathematics Subject Classification}
\makeatother

\title[Filtered Boolean powers of simple algebras]{Filtered Boolean powers of finite simple non-abelian Mal'cev algebras}
\author{Peter Mayr}
\address{Department of Mathematics,
University of Colorado Boulder, USA}
\email{peter.mayr@colorado.edu}
\author{Nik Ru{\v s}kuc}
\address{School of Mathematics and Statistics,
University of St Andrews, Scotland, UK}  
\email{nik.ruskuc@st-andrews.ac.uk}
\thanks{Supported by the Engineering and Physical Sciences Research Council EP/V003224/1, by the Heilbronn Institute for Mathematical Research, and by the National Science Foundation under Grant No. DMS 2452289.}


\date{\today}

\begin{document}
\maketitle

\begin{abstract}
 Let $\A$ be a finite simple non-abelian Mal'cev algebra (e.g.\ a group, loop, ring).
 We investigate the Boolean power $\D$ of $\A$ by the countable atomless Boolean algebra $\B$ filtered at 
 some idempotents $e_1,\dots,e_n$ of $\A$.
When $e_1,\dots,e_n$ are \emph{all} idempotents of $\A$ we establish two concrete representations of $\D$: as 
the (generalized) Fra{\"i}ss{\'e} limit of the class of finite direct powers of $\A$, and as
congruence classes of the
 countable free algebra in the variety generated by $\A$.
 Further, for arbitrary $e_1,\dots,e_n$, we show that $\D$ is $\omega$-categorical and that its automorphism group
 has the small index property, strong uncountable cofinality and the Bergman property.
 As necessary background we establish some general properties of congruences and automorphisms of filtered Boolean powers of $\A$ by any Boolean algebra $\B$, including a semidirect decomposition for their automorphism groups.
\end{abstract}

\section{Introduction and main results}

The purpose of this paper is to investigate filtered Boolean powers~$\D$ of a finite simple non-abelian Mal'cev algebra $\A$ by a Boolean algebra~$\B$.
On a general level, we establish some properties of congruences and automorphisms of $\D$.
When $\B$ is the countable atomless Boolean algebra, we investigate the role $\D$ plays in the variety generated by $\A$, and prove that it has a number of combinatorial properties: 
it is $\omega$-categorical and its automorphism group has the small index property, strong uncountable cofinality and the Bergman property.

In this section we present the background and motivation for our work, introduce the concepts and notation that will be used throughout and give statements of the main results.

 A \emph{variety} is a class of algebraic structures (\emph{algebras} for short) of the same type that is defined by
 equations. By Birkhoff's Theorem the \emph{variety generated by a class of algebras} $C$ is the class of all
 homomorphic images of subalgebras of products of elements in $C$. We refer to~\cite{BS:CUA} for background on
 general algebra.
  
 Varieties generated by a finite simple group have been studied by Neumann~\cite{Ne:VG}, Apps~\cite{Ap:ACG} and others, 
 while Apps \cite{Ap:BPG} considers subgroups of Boolean products and powers.
 Many of the techniques carry over from groups to Mal'cev algebras in general. Here an algebra $\mathbf{A}$ is
 \emph{Mal'cev} if it has a ternary term operation $m$ satisfying $m(x,x,y) = y = m(y,x,x)$.
 Note that groups, quasigroups, loops, rings as well as their expansions
with additional operations
 are all Mal'cev.

 The generalization of commutator theory from groups to algebras as described in~\cite{FM:CTC} allows us to talk
 about (non)-abelian algebras. A Mal'cev algebra is \emph{abelian} if its basic operations can be represented
 as affine functions over a module; otherwise we call it \emph{non-abelian}. This is consistent with the classical
 notions for groups and loops. A ring is abelian if and only if its multiplication is constant $0$.

 Unless otherwise indicated, $\A$ will always stand for a finite simple non-abelian Mal'cev algebra.
 Throughout, we will use bold letters to denote algebras and other structures,
 and matching standard letters will be used for the underlying universes. Thus, here, $A$ would be the universe of $\A$.
 Let $V$ denote the variety generated by $\A$, and let $W$ be the variety generated by all proper subalgebras of $\A$.
 Then $W$ is the unique maximal subvariety of $V$.
 Every finite algebra in $V$ is isomorphic to $\A^k\times\B$ for
 some integer $k\geq 0$
 and some finite $\B\in W$ (cf.~\cite{Ne:VG} for groups, \cite{Fo:GBT2} for functionally complete algebras).
 This indicates the distinguished role the class $\{ \A^k\st k\geq 1\}$ of finite powers of $\A$ plays in the
 structure of the finite members of $V$.
 It turns out that this class has the amalgamation property, leading to the existence of a (generalized)
 Fra{\"i}ss{\'e} limit.
 
In order to be able to state this, we briefly review the standard Fra{\"i}ss{\'e} theory following \cite{Ho:MT}.
Let $K$ be a countable class of finitely generated structures of the same type.
Consider the following three properties that such a class may have:
\vspace{1mm}

\noindent
\textit{Hereditary property (HP):} Every finitely generated substructure of a member of $K$ is isomorphic to a member of $K$.
\vspace{1mm}

\noindent
\textit{Joint embedding property (JEP):} For any $\A,\B\in K$ there exists $\C\in K$ into which both $\A$ and $\B$ embed.
\vspace{1mm}

\noindent
\textit{Amalgamation property (AP):} For any $\A,\B,\C\in K$ with embeddings $\varphi\colon\A\to\B$ and $\psi\colon \A\to\C$ there exists $\D\in K$ and 
embeddings $\varphi'\colon\B\to \D$ and $\psi'\colon\C\to\D$ such that $\varphi'\varphi=\psi'\psi$.
\vspace{1mm}

The classical Fra{\"i}ss{\'e} Theorem asserts that these three conditions are equivalent to the existence of
 a countable homogeneous structure whose class of finitely generated substructures is $K$. This necessarily unique
 structure is known as the \emph{Fra{\"i}ss{\'e} limit} of $K$; see \cite[Theorem 7.1.2]{Ho:MT}.

 Returning to our situation of the variety $V$ generated by
 a finite simple non-abelian Mal'cev algebra $\A$, note that its subclass of finite algebras
 has the HP but in general not the JEP or AP.
On the other hand the subclass $K=\{ \A^k\st k\geq 1\}$ 
will be shown to have the JEP and AP but in general not the HP.
Under these conditions,  a Fra{\"i}ss{\'e} limit-like structure still exists with similar properties as the classical
Fra{\"i}ss{\'e} limit.
We call this structure the \emph{generalized Fra{\"i}ss{\'e} limit} of $K$.
More details are given in Subsection \ref{subsec:GenFre}.

This generalized Fra{\"i}ss{\'e} limit can be explicitly described as a filtered Boolean power of $\A$,
as used by Arens and Kaplansky~\cite{AK:TRA} (see also~\cite{Ev:EAC}).
 We briefly review this construction.

Let $\B$ be a Boolean algebra. The \emph{Stone space} of $\B$ is the set $X$ of ultrafilters on $\B$
 with the topology whose basic open sets are $\{ x\in X\st b\in x\}$ for $b\in B$.
By Stone's Representation Theorem $\B$ is isomorphic to the Boolean algebra of clopen sets of $X$.
An arbitrary topological space $X$ is homeomorphic to the Stone space of a Boolean algebra if and only if it is compact and totally separated.
For basics on Boolean algebras we refer the reader to \cite[Chapter 1]{Mo:HBA1} and for Stone duality to \cite[Chapter 3]{Mo:HBA1}.

Recall that there is a unique (up to isomorphism) countable atomless Boolean algebra $\B$.
 It is in fact the  Fra{\"i}ss{\'e} limit of all finite Boolean algebras, and will play a particularly important role in this paper.
 Its Stone space is the \emph{Cantor space} $2^\omega$,
 i.e. the unique space up to homeomorphism which is compact, Hausdorff, has no isolated points and has a countable basis of clopen sets.

Let $\B$ be a Boolean algebra with Stone space $X$, and let $\A$ be an algebra. 
The \emph{Boolean power} $\A^\B$ is the subalgebra of $\A^X$ with universe
 \[ A^\B := \{ f\colon X\to A \st f \text{ is continuous} \}, \]
 where $A$ is endowed with the discrete topology.

 Borrowing the concept from semigroups, we call an element $e$ of an algebra $\A$ an \emph{idempotent}
 if $\{e\}$ forms a (singleton) subalgebra of $\A$.  
 For pairwise distinct $x_1,\dots,x_n \in X$ and idempotents $e_1,\dots,e_n$ of $\A$,
the subalgebra $\ABxe$ of $\A^\B$ with universe
\[ \{ f\in A^\B \colon f(x_1)=e_1,\dots,f(x_n)=e_n \} \]
 is a \emph{filtered Boolean power} of $\A$. 
 
 In choosing our notation we have opted to follow the usage in group theory literature, for example \cite{Ap:BPG,BE:SIP,BG:ARFG}, because it emphasizes the power nature of the construction. It is, however, important not to confuse this construction with the standard direct power $\A^Y$ of $|Y|$ copies of $\A$, where $Y$ is an arbitrary (unstructured) set.
 Other notations used for $\A^\B$ are e.g.  $\A[\B]$ in \cite{Bu:BP,Ev:EAC} or $\A[\B]*$ in \cite{BS:CUA}.
 We can now state our first main result.

\begin{thm} \label{thm:Flim}
 Let $\A$ be a finite simple non-abelian Mal'cev 
 algebra, let the idempotents of $\A$ be $e_1,\dots,e_n$,
let $\B$ be the countable atomless Boolean algebra,
and let $x_1,\dots,x_n\in 2^\omega$ be distinct.
\begin{enumerate}
\item \label{it:Flim}
The class $\{\A^k\st k\geq 1\}$ has the JEP and AP and its generalized Fra{\"i}ss{\'e} limit is isomorphic to $\ABxe$.
\item \label{it:free}
 Let $\F$ be the free algebra of countable rank in the variety $V$ generated by $\A$, and let $\theta$ be the smallest
 congruence such that $\F/\theta$ is in the variety $W$ generated by all proper subalgebras of~$\A$.
 Then $\F$ embeds into $\A^\B$ and every $\theta$-class which is a subalgebra of $\F$ is isomorphic to
 $\ABxe$.
\end{enumerate}  
\end{thm}

 In the theorem above it is in fact sufficient to take $e_1,\dots,e_n$ to be a set of orbit representatives under $\Aut\A$ for
 the idempotents in $\A$ (see Corollary~\ref{cor:isoreduced}).
 
 Part \eqref{it:Flim}  of Theorem~\ref{thm:Flim} is proved in Section~\ref{sec:Flim} and 
 part \eqref{it:free} in Section~\ref{sec:variety}.

Note that when $\A$ is a group, the filtered Boolean power arising in Theorem~\ref{thm:Flim}\eqref{it:free} is of the form
$(\A^\B)^x_e$, where $e$ is the identity of~$\A$. This is in turn easily seen to be isomorphic to $\A^\R$ where $\R$ is the countable atomless Boolean ring without identity. (For Boolean powers over Boolean rings see \cite{Ap:BPG,HiPl:BPDP}.) 
Theorem~\ref{thm:Flim}\eqref{it:free}
in this case gives that the kernel of the homomorphism from $\F$ to the free group of countable rank in $W$ is isomorphic to $(\A^\B)^x_e$. This observation was made in \cite[p.\ 367-8]{BE:SIP} and \cite[p.\ 201]{BG:ARFG}.
As a consequence of Theorem~\ref{thm:Flim}\eqref{it:free},
for groups and  loops and rings we obtain
the stronger result that every countable algebra $\C$ in $V$ is a split
 extension of a Boolean power of $\A$ by an algebra in $W$ (Corollary~\ref{cor:group}).

 We also provide an illustration of Theorem~\ref{thm:Flim} for a specific algebra with more than one idempotent.

\begin{exa} \label{exa:twoids}
 Let $\A := (\Z_2,x-y+z,\cdot)$ be the idempotent reduct of the field of size $2$.
 The only proper subalgebras of $\A$ are the trivial ones $\{0\}$ and $\{1\}$, and $\Aut\A$ is trivial as well.
 It turns out that the variety $V$ generated by $\A$ is the class of (isomorphic copies of) filtered Boolean
 powers of $\A$.

 Every finite algebra in $V$ is isomorphic to $\A^k$ for some integer $k\geq 0$. 
 Let $\B$ be the countable atomless Boolean algebra. Then the generalized Fra{\"i}ss{\'e} limit of $\{ \A^k \st k\geq 1 \}$
 is isomorphic to the free algebra of countable rank in $V$ and to
\[ (\AB)_{0,1}^{x_0,x_1} \cong (\AB)_0^{x_0} \times (\AB)_1^{x_1} \]
 for distinct $x_0,x_1\in 2^\omega$.
\end{exa} 

 We now investigate automorphism groups of filtered Boolean powers.

 It is well-known that for any countably infinite first order structure~$M$, 
 the group of automorphisms
 $G := \Aut M$ affords a natural \emph{topology of pointwise convergence} with basic open sets the cosets of
 pointwise  stabilizers
$G_Y$ of finite subsets $Y$ of $M$.
 This makes $G$ a topological (in fact, Polish) group.

 The group $G$ has the \emph{small index property} (SIP) if any subgroup of index less than
 $2^{\aleph_0}$ is open, that is, contains the pointwise stabilizer of some finite subset of $M$.

 Note that if $\Aut M$ has the SIP, then its topological structure is completely determined by its abstract
 algebraic structure.
Specifically, if $N$ is another countable structure and $\Aut M\cong\Aut N$ as groups, then in fact $\Aut M$ and $\Aut N$ are isomorphic as topological groups ~\cite[Proposition 5.2.2]{Mac:SHS}.
 In particular $\Aut M$ has a unique Polish topology.

 The \emph{cofinality} of a group $G$ is the smallest cardinality $\kappa$ such that $G$ is the union of a chain of length
 $\kappa$ of proper subgroups. 

 The \emph{strong cofinality} of $G$ is the smallest cardinality $\kappa$ such that $G$ is the union of a chain
 $(U_i)_{i<\kappa}$ of proper subsets of $G$ such that for all $i<j$ we have $U_i\subseteq U_j$, 
 $U_i = U_i^{-1}$ and $U_iU_i \subseteq U_k$ for some $k\in\kappa$. 

 The group $G$ has the \emph{Bergman property} if for each subset $E$ of $G$ such that $1\in E = E^{-1}$ and
 $E$ generates $G$ there exists $k\in\N$ such that $E^k = G$.

 Finally $G$ has \emph{ample generics} if for any $k\in\N$ the diagonal conjugation action of $G$ on $G^k$
 has a comeager orbit, that is, an orbit containing the intersection of countably many dense open subsets of $G^k$.

 These properties are related as follows.
 A group $G$ has uncountable strong cofinality if and only if 
 $G$ has uncountable cofinality and the Bergman property
 (cf.~\cite[Theorem 2.2]{DH:GAG}).
 The existence of ample generics implies the SIP as shown by Hodges, Hodkinson, Lascar and Shelah~\cite{HHLS:SIP}.
 Further, if $M$ is $\omega$-categorical and $\Aut M$ has ample generics,
 then $\Aut M$ has uncountable strong cofinality by Kechris and Rosendal~\cite{KR:TAG}.

 The automorphism group $G$ of the countable atomless Boolean algebra
(equivalently, the group of homeomorphisms of the Cantor space~$2^\omega$)
 was among the first for which these
 properties were established. Truss \cite{Tr:IPG2} showed that it has the SIP. Droste and G\"obel~\cite{DG:UCPG}
 proved uncountable strong cofinality.
 Later Kwiatkowska~\cite{Kw:GHC} subsumed these results by showing that $G$ actually has ample generics.

 As our second main result we show that the SIP and uncountable strong cofinality carry over from the countable atomless Boolean algebra $\B$ to filtered Boolean powers of a finite simple non-abelian Mal'cev algebra by $\B$.
   
\begin{thm} \label{thm:SIP}
 Let $\A$ be a finite simple non-abelian Mal'cev algebra, let $n\geq 0$, let $e_1,\dots,e_n\in A$ be idempotents,
 let $\B$ be the countable atomless Boolean algebra, and let $x_1,\dots,x_n\in 2^\omega$ be distinct. Then 
\begin{enumerate} 
\item \label{it:omegacat}
  $\ABxe$ is $\omega$-categorical;
\item \label{it:SIP}
$\Aut\ABxe$ has the SIP;
\item \label{it:Bergman}
 $\Aut\ABxe$ has uncountable strong cofinality, in particular, uncountable cofinality and the Bergman property.
\end{enumerate}  
\end{thm}

Theorem~\ref{thm:SIP}\eqref{it:omegacat} will follow from results of Macintyre and Rosenstein~\cite{MR:CRN}
 in Section~\ref{sec:omegacat}.
 We will prove items~\eqref{it:SIP},\eqref{it:Bergman} of Theorem~\ref{thm:SIP} in
 Sections~\ref{sec:SIP},~\ref{sec:Bergman} following the strategies for the countable atomless Boolean
 algebra by Truss~\cite{Tr:IPG2} and Droste and G\"obel~\cite{DG:UCPG}, respectively.
 Underpinning these results is an explicit description of $\Aut\ABxe$ as semidirect product of
 the pointwise stabilizer of $x_1,\dots,x_n$ in the group of homeomorphisms of $2^\omega$
and a normal subgroup $K$, which is a certain
 closure of a filtered Boolean power of $\Aut\A$ by $\B$
 (see Theorem~\ref{thm:AutABxe} and the subsequent remarks).
 
 It remains open whether $\Aut\ABxe$ has ample generics. In particular, it is unclear whether Kwiatkowska's
 approach for the countable atomless Boolean algebra from~\cite{Kw:GHC} can be extended to filtered Boolean powers.\footnote{While this paper was under review, the authors showed in \cite{MR:AGBP} that $\Aut\ABxe$ as in Theorem~\ref{thm:SIP} has ample generics based on Theorem~\ref{thm:AutABxe} and a generalization of Kwiatkowska's projective Fra{\"i}ss{\'e} limit construction for generics in the homeomorphism group of $2^\omega$ \cite{Kw:GHC}.}

 \section{Preliminaries} \label{sec:prelim}

 In this section we gather general terminology and notation that is needed throughout the paper. More specialised concepts and results will be introduced where they are needed in the text.

  \subsection{Notation}
  The set of natural numbers $\{1,2,\dots\}$ will be denoted by $\N$.
 For $n\in\N$ we write $[n] := \{1,\dots,n\}$.
First order structures will be denoted by the boldface capital letters $\A$, $\B$, $\C$ etc.,
with the corresponding letters $A$, $B$, $C$ standing for their underlying sets (universes). One exception will be groups, which will be denoted by non-bold letters
$G$, $H$ etc. 
For an algebra $\A$ and a subset $U\subseteq A$ we will denote by $\langle U\rangle$ the \emph{subuniverse} of $\A$ \emph{generated} by $U$, i.e. the carrier set of the subalgebra of $\A$ generated by $U$.
Topological spaces will be denoted by $X$, $Y$, $Z$ etc. For a structure $\A$ we will write $\Aut\A$ for its automorphism group, and for a topological space $X$ we will write $\Homeo X$ for its group of homeomorphisms. For a permutation group acting on a set $X$ and a subset $Y\subseteq X$ the pointwise stabilizer of $Y$ in $G$ is denoted by
$G_Y  := \{ g\in G \st g(x)=x \text{ for all } x\in Y \}$.

\subsection{Boolean algebras}

Stone's Representation Theorem will be used throughout the paper, usually without explicit reference. An abstract Boolean algebra $\B$ will be at times identified as the algebra of clopen sets of its Stone space $X$.
It will be useful to bear in mind that $\Aut\B$ and $\Homeo X$ are naturally isomorphic to each other; see \cite[Theorem 8.2]{Mo:HBA1}.

\subsection{Boolean powers}

Our main references are \cite{Bu:BP} and \cite[Section IV.5]{BS:CUA}.
Our definition follows the latter source.
In \cite{Bu:BP} instead of continuous functions from the Stone space $X$ of $\B$ into $\A$, the elements are taken to be
mappings from $A$ to $B$. To obtain such a mapping from a continuous function $f\colon X\rightarrow A$, one just maps
$a\in A$ to its pre-image $f^{-1}(a)$, regarded as an element of $B$ by Stone's duality. Equipped with this it is easy to
see that the two definitions are equivalent.

 Regarding \emph{filtered} Boolean powers, there is no universally agreed terminology or notation in literature. 
 The concept first seems to appear in~\cite{AK:TRA} where Arens and Kaplansky use it for representing rings.
 In~\cite{Fo:GBT1} Foster generalizes it to arbitrary algebras calling it Boolean extensions.
 As far as we know the first use of the term filtered Boolean powers is due to Burris and Werner~\cite{BW:SCT}. 
 This term is also preferred by Evans in his survey \cite{Ev:EAC}.

\subsection{Subdirect products of Mal'cev algebras} 
 Certain structural features of subgroups of direct products of groups generalize to Mal'cev algebras.
 We present those that will be used in this paper. 
 First is Fleischer's extension of Goursat's Lemma for groups~\cite[Theorem 5.5.1]{Ha:TG}.

\begin{lem}[cf. Fleischer's Lemma {\protect{\cite[Lemma IV.10.1]{BS:CUA}}}]
\label{la:Fleischer}
Let $\A_1,\A_2$ be algebras belonging to a variety $V$ with a Mal'cev term, and
let $\C\leq\A_1\times\A_2$ be a subdirect product of $\A_1$ and $\A_2$.
Then there exist $\D\in V$ and epimorphisms $\psi_i\colon \A_i\rightarrow\D$ for $i\in [2]$ such that
\[
C=\{ (a_1,a_2)\in A_1\times A_2\colon \psi_1(a_1)=\psi_2(a_2)\}.
\]
\end{lem}




The second states that a finite subdirect power of a simple Mal'cev algebra has a very restricted form:
it is entirely determined by the projections to a subset of direct components, and these projections are pairwise independent.

\begin{lem}[cf. Foster--Pixley Theorem~{\protect{\cite[Corollary IV 10.2]{BS:CUA}}}]
\label{la:FP}
Let $\A$ be a simple Mal'cev algebra, let $n\in\N$, and let $\C\leq \A^n$ be a subdirect power.
Then there exist a non-empty subset $I\subseteq [n]$, a mapping $f\colon {[n]\setminus I}\rightarrow I$ and automorphisms $\theta_i\in\Aut\A$ for $i\in [n]\setminus I$ such that
\[
C=\{ (a_1,\dots,a_n)\in A^n\colon a_i=\theta_i(a_{f(i)})\text{ for all }i\in [n]\setminus I\}.
\]
In particular, $\C\cong \A^{|I|}$.
\end{lem}

Finally, a direct power of a simple non-abelian Mal'cev algebra has no skew-congruences.

\begin{lem}[cf.~{\protect\cite[Lemma 9.69, Corollary 9.66]{FMMT:ALV3}}]
\label{la:noskew}
Let $\A$ be a finite simple non-abelian Mal'cev algebra and let $n\in\N$.
Every congruence of the direct power $\A^n$ is a product congruence,
i.e. it has the form $\theta_1\times\dots\times \theta_n$, where each $\theta_i$ is either the total relation or the equality on~$\A$.
\end{lem}

\subsection{Generalized Fra{\"i}ss{\'e} limits.}
\label{subsec:GenFre}

As mentioned in the introduction, some of our filtered Boolean powers will be representable as so called generalized Fra{\"i}ss{\'e} limits of finite powers.
The need for generalization comes from the fact that the finite powers do not satisfy the hereditary property from the classical Fra{\"i}ss{\'e} theory.

\begin{lem}[{\protect\cite[Exercise 7.1.11]{Ho:MT}}]
\label{la:genFr}
Let $K$ be a countable class of first order structures of the same type satisfying the JEP and AP.
Then there exists a countable structure $\D$ which satisfies the following properties: 
\begin{enumerate}
\item\label{it:gfr1}
up to isomorphism, the class of all finitely generated substructures of $\D$ is the class of all finitely generated substructures of members of $K$;
 \item\label{it:gfr2}
 $\mathbf{D}$ is a direct limit of structures in $K$; and
 \item\label{it:gfr3}
 every isomorphism between substructures of $\mathbf{D}$ that are isomorphic to some element in $K$ extends
 to an automorphism of $\mathbf{D}$.
 \end{enumerate}
 \end{lem}
 
 We will call the structure $\D$ in the above result the \emph{generalized Fra{\"i}ss{\'e} limit} of $K$.
 Property \eqref{it:gfr3} will be referred to as \emph{$K$-homogeneity}.
 

\section{Filtered Boolean powers} \label{sec:Booleanpowers}

We collect the structural information on filtered Boolean powers and their automorphism groups that we need for
proving our main results.
 Throughout this section we use the following notation unless specified otherwise:
\begin{itemize}
\item $\A$ is a finite simple non-abelian Mal'cev algebra;
\item $n\in\N\cup\{0\}$;
\item $e_1,\dots,e_n$ are idempotents of $\A$; they do not need to be distinct or to include all the idempotents;
\item
  $\B$ is a Boolean algebra with Stone space $X$; the latter may be taken to consist of the ultrafilters on $\B$;
 \item
$x_1,\dots,x_n$ are pairwise distinct points in $X$, and 
 $X^\circ := X\setminus\{x_1,\dots,x_n\}$ is equipped with the subspace topology from $X$.
\end{itemize}

\subsection{Filtered Boolean powers represented on arbitrary sets}

 Filtered Boolean powers can be defined using an arbitrary representation of a Boolean algebra as field of sets,
 not only via its Stone space.

\begin{lem} \label{lem:fBp}
 Let $\B$ be a subalgebra of the Boolean algebra of subsets of a set $Y$. 
\begin{enumerate}
\item \label{it:y'} For any $y\in Y$ the set $y' := \{b\in B \st y\in b \}$ is an ultrafilter on $\B$.
\item \label{it:fBp} Let $\A$ be a finite algebra, let $e_1,\dots,e_n\in A$ be idempotents, and let $y_1,\dots,y_n\in Y$
 be distinct. Then the set
  \begin{align*}
D^{y_1,\dots,y_n}_{e_1,\dots,e_n} & :=
  \{ f\in A^Y \st f^{-1}(a)\in B \text{ for all } a\in A,\\
 &\hspace{30mm}\text{and } f(y_i) = e_i \text{ for all } i\in [n] \}
  \end{align*}
 \medskip
\noindent
 is the universe of a subalgebra $\D^{y_1,\dots,y_n}_{e_1,\dots,e_n}$ of $\A^Y$
 isomorphic to the filtered Boolean power $(\A^\B)^{y'_1,\dots,y'_n}_{e_1,\dots,e_n}$.
\end{enumerate}
\end{lem}  

\begin{proof}
 \eqref{it:y'} is clear.
  
 \eqref{it:fBp} For $b\in B$, let $b^{**} := \{ x\in X \st b\in x \}$. Let $B^{**} := \{ c\subseteq X \st c \text{ clopen} \}$.
 By Stone's Representation Theorem
 \[ \varphi\colon \B\to\B^{**},\ b \mapsto b^{**}, \]
 is an isomorphism. Now
\[ D := \{ f\in A^Y \st f^{-1}(a)\in B \text{ for all } a\in A \} \]
is the universe of a subalgebra $\D$ of $\A^Y$.
 For $f\in D$ define $f'\colon X \to A$ by
\[ (f')^{-1}(a) = (f^{-1}(a))^{**} \text{ for all } a\in A. \]
 Equivalently, if $Y = b_1\sqcup \dots\sqcup b_k$ (disjoint union) and $f(b_j) = \{a_j\}$ for $b_j\in B, j\in [k]$, then
 $X = b^{**}_1\sqcup \dots\sqcup b^{**}_k$ and $f'(b^{**}_j) = \{a_j\}$ for $j\in [k]$. 
 Since $\varphi\colon b\mapsto b^{**}$ is an isomorphism and $f(b)=f'(b^{**})$ for all $b\in B$,
 it follows that
\[ \varphi'\colon \D \to \A^\B,\ f\mapsto f', \]
is an isomorphism.
 
  We claim that $\varphi'$ restricts to the required isomorphism between
 $\D^{y_1,\dots,y_n}_{e_1,\dots,e_n}$ and $(\A^\B)^{y'_1,\dots,y'_n}_{e_1,\dots,e_n}$. To see that 
 the domain and codomain are correct,
 note that for $f \in D$ and $i\in [n]$,
\begin{align*} f(y_i) = e_i &\Leftrightarrow y_i\in f^{-1}(e_i) \Leftrightarrow  f^{-1}(e_i) \in y'_i\\
& \Leftrightarrow y'_i\in (f^{-1}(e_i))^{**} = (f')^{-1}(e_i)
 \Leftrightarrow f'(y_i') = e_i.
\end{align*}
 Hence 
 \[ \varphi'(\D^{y_1,\dots,y_n}_{e_1,\dots,e_n})  =  (\A^\B)^{y'_1,\dots,y'_n}_{e_1,\dots,e_n}, \]
 completing the proof.
\end{proof}

\subsection{Congruences}
 Our description of congruences of filtered Boolean powers extends those of Boolean powers of simple algebras by
 Burris~\cite[Theorem 3.5, Corollary 3.6]{Bu:BP}.

For $\D := \ABxe$ and $Y\subseteq X$ let
\[ \theta_Y := \{ (f,g)\in D^2 \st f|_{Y} = g|_{Y} \} \]
 denote the kernel of the projection of $\D$ to $Y$.

\begin{lem}[{cf.\ {\cite[Theorem 3.5]{Bu:BP}}}] \label{lem:con}
 Let $\D  := \ABxe$.
\begin{enumerate}  
\item \label{it:principal}
  The principal congruence of $\D$ generated by a pair $(f,g)\in D^2$ is the kernel
  of the projection of $\D$ onto
 the clopen set
\[ b:=\{ x\in X \st f(x)=g(x) \}. \]
\item \label{it:con} 
Every congruence $\theta$ of $\D$ is the kernel of the projection of $\D$ onto
 \[ Y := \bigcap_{(f,g)\in\theta} \{ x\in X \st f(x)=g(x) \}. \]
\end{enumerate} 
\end{lem}

\begin{proof}
\eqref{it:principal}
 Denote by $\theta$ the principal congruence generated by $(f,g)$.
 Note that every finite subalgebra of $\D$ is contained in a filtered Boolean power $\D\cap \A^\C$ for some finite subalgebra
 $\C$ of $\B$. Here we view $C$ as a finite set of clopens on $X$ and $\A^\C$ as the subalgebra of $\A^\B$ with the universe
\[ A^\C = \{h\colon X\to A \st h^{-1}(a) \in C \text{ for all } a\in A\}. \]

Now let $\C$ be a finite subalgebra of $\B$ such that $\A^\C$ contains both $f$ and $g$.
 For all $x\in b$ the projection of $\theta\cap (D\cap A^\C)^2$ onto $x$ is contained in the equality on $A$;
 else if $x\in X\setminus b$ we have $f(x)\neq g(x)$, and hence the projection of $\theta\cap (D\cap A^\C)^2$ onto $x$
 is the total congruence on $\A$, as $\A$ is simple.
 By Lemma \ref{la:noskew} each congruence of the finite algebra $\D\cap \A^\C$ is a product congruence
 and uniquely determined by its projections onto the points $x\in X$. It follows that
 \[  \theta \cap (D\cap A^\C)^2  = \{ (u,v)\in (D\cap A^\C)^2 \st u|_b  = v|_b \}. \]
 Hence the restrictions of $\theta$ and the projection kernel $\theta_{b}$ onto any finite $\D\cap\A^\C$
 containing $f$ and $g$ coincide.
 This proves~\eqref{it:principal}.

\eqref{it:con}  
 For $S$ a set of clopens in $X$, we claim that 
 \begin{equation}  \label{eq:thetaS}
   \bigvee_{b\in S} \theta_b = \theta_{\bigcap S}.
 \end{equation} 
 From Lemma \ref{la:noskew} it is straightforward that $\theta_b \vee \theta_c = \theta_{b\cap c}$ for clopens $b,c$ on $X$.
 Hence $(f,g)\in \bigvee_{b\in S} \theta_b$ if and only if $f$ and $g$ agree on some element in the filter $F$ of clopens
 on $X$ that is generated by $S$.
 As a standard consequence of the Ultrafilter Lemma, every filter is the intersection of the ultrafilters which contain it.
 By Stone's Representation Theorem and Lemma~\ref{lem:fBp}\eqref{it:y'},
 the ultrafilters of the Boolean algebra of clopens on $X$ are exactly of the form
 $\{ b \text{ clopen} \st x\in b \} =: x'$ for $x\in X$.
 Then $F \subseteq x'$ if and only if every $b\in F$ contains $x$ if and only if $x\in\bigcap S$.
 Hence
\[
 F = \bigcap \{ x' \st x\in\bigcap S \} = \{ b \text{ clopen } \st \bigcap S \subseteq b \}
\]
and \eqref{eq:thetaS} follows. Since
\[ \theta = \bigvee_{(f,g)\in\theta} \theta_{\{ x\in X \st f(x)=g(x) \}} \]
using~\eqref{it:principal}, we obtain $\theta = \theta_Y$ as claimed.
\end{proof}

 For a subset $Y$ of $X$, the quotient $\ABxe/\theta_Y$ is isomorphic to the restriction of $\ABxe$ to $Y$.
 The latter is again isomorphic to a filtered Boolean power by Lemma~\ref{lem:fBp}.
 Thus every homomorphic image of a filtered Boolean power is a filtered Boolean power again.

 Lemma~\ref{lem:con} actually yields a lattice isomorphism between the lattice of congruences
 of the filtered Boolean power $\ABxe$ and the lattice of filters contained in the intersection of ultrafilters
 $\bigcap_{i=1}^n x_i$ on the Boolean algebra $\B$. In particular every principal filter generated by some clopen
 $b$ containing the points $x_1,\dots,x_n$ corresponds to a principal congruence on $\ABxe$ generated by a pair of
 functions $f,g$ that coincide exactly on $b$.

\subsection{Automorphisms} \label{sec:autos}

In this subsection we give a detailed analysis of the automorphisms of $\D := \ABxe$,
in terms of automorphisms of $\A$ and of $\B$.
With some additional mild assumptions regarding $e_1,\dots,e_n$ we obtain a semidirect decomposition of 
$\Aut \D$ in Theorem~\ref{thm:AutABxe}.
 Our results generalize those for Boolean powers of groups by Boolean rings
 by Apps~\cite[Theorem C]{Ap:BPG}.

 Automorphisms of $\B$ act naturally on the ultrafilters of $\B$, 
 i.e.\ the points in  the Stone space $X$ of $\B$.
This correspondence witnesses the well-known fact that $\Aut\B\cong\Homeo X$.
 Reversing it and lifting to the Boolean power leads to the following construction.
 For $\psi\in (\Homeo X)_{\{x_1,\dots,x_n\}}$ define
\begin{equation} \label{eq:defpsiD}
 \psi^\D\colon D\to D,\ f\mapsto f\psi^{-1}. 
\end{equation}
Equivalently $[\psi^\D(f)]^{-1} (a) = \psi f^{-1}(a)$ for all $f\in D, a\in A$.
 It is easy to check that $\psi^\D$ is an automorphism of $\D$ and that
\begin{equation} \label{eq:defg}
  g\colon (\Homeo X)_{\{x_1,\dots,x_n\}} \to \Aut\D,\ \psi\mapsto\psi^\D,
\end{equation}
is a group embedding.

In the following lemma we show that, conversely, automorphisms of~$\D$ induce homeomorphisms on $X^\circ$, using the fact that every automorphism of $\D$ induces an automorphism on the congruence lattice of~$\D$.

\begin{lem} \label{lem:hphi}
 For $x\in X$, let $\pi_x\colon \ABxe\to \A,\ f\mapsto f(x)$, be the projection at $x$,
 and let $\theta_x:=\ker\pi_x$ be its kernel. For any $\varphi\in\Aut\ABxe$ the following hold:
\begin{enumerate}
\item \label{it:phicirc1}
 The mapping $\varphi^{\circ}\colon X^\circ\to X^\circ$ defined by
\[ \varphi(\theta_x) = \theta_{\varphi^\circ(x)} \]
 is a homeomorphism.
\item \label{it:phix}
 For every $x\in X^\circ$, the mapping $\varphi_x\colon A\rightarrow A$ defined by
\begin{equation}
\label{eq:phixdef}
\varphi_x \pi_x = \pi_{\varphi^\circ(x)} \varphi
\end{equation}
is an automorphism of $\A$.
\item \label{it:phicirc2}
 For any convergent net $(y_i)_{i\in I}$ in $X^\circ$ with limit $x_k$ for $k\in[n]$, all cluster points of
 $(\varphi^\circ(y_i))_{i\in I}$ are in $\{x_\ell \st \ell\in [n],\ e_\ell\in(\Aut\A)(e_k)\}$.
\end{enumerate}  
\end{lem}

\begin{proof}
 Let $\D := \ABxe$.
  
\eqref{it:phicirc1}
 To see that $\varphi^\circ$ is well-defined, let $x\in X^\circ$.
 Since  $\A$ is simple, $\theta_x$ is a maximal congruence of $\D$, which $\varphi$ maps to another maximal congruence
\[ \varphi(\theta_x) = \{ (\varphi(f),\varphi(g))\in D^2 \st f(x)=g(x) \}. \]
 By Lemma~\ref{lem:con} we have a unique point $\varphi^\circ(x)$ in $X^\circ$ such that
 $\varphi(\theta_x) = \theta_{\varphi^\circ(x)}$ as required.
 
 Next note that $\varphi^\circ$ is bijective with inverse $(\varphi^{-1})^\circ$ since by
 definition
\[ \theta_x = \varphi^{-1}(\theta_{\varphi^\circ(x)}). \]

 For continuity, it suffices to show that $\varphi^\circ$ maps any clopen subset $b$ of $X$ that is contained in $X^\circ$
 to such a clopen subset again.
 As mentioned in the comment after Lemma~\ref{lem:con}, the kernel of the projection on the clopen
  $X\setminus b$ is the principal congruence generated by a pair of functions $f,g\in D$ which coincide exactly on
  $X\setminus b$.
It follows that the congruence generated by $(\varphi(f),\varphi(g))$ is the kernel of the projection to the clopen
\[ \{x\in X \st \varphi(f)(x)=\varphi(g)(x) \} = X\setminus\varphi^\circ(b). \]
 Hence $\varphi^\circ(b)$ is clopen in $X$ and disjoint from $\{x_1,\dots,x_n\}$.

 \eqref{it:phix}
 Let $x\in X^\circ$.
Note that \eqref{eq:phixdef} can be rewritten as
\[
\varphi_x(f(x)) = \varphi(f)(\varphi^\circ(x))\quad\text{for all } f\in D.
\] 
 Suppose $f,g\in D$ satisfy $f(x) = g(x)$, that is, $(f,g)\in\theta_x$.
 Then $(\varphi(f),\varphi(g))$ is in $\varphi(\theta_x) = \theta_{\varphi^\circ(x)}$.
 Hence
\[ \varphi_x(f(x)) = \varphi(f)(\varphi^\circ(x)) = \varphi(g)(\varphi^\circ(x)) = \varphi_x(g(x)). \]
 So~\eqref{eq:phixdef} uniquely defines $\varphi_x$.

To check that $\varphi_x$ is a homomorphism on $\A$, let $t$ be a $k$-ary operation in the signature of $\A$ and
$f_1,\dots,f_k\in D$. Then
\begin{align*}
  \varphi_x(t^\A(f_1(x),\dots,f_k(x))) & = \varphi(t^\D(f_1,\dots,f_k))\, (\varphi^\circ(x)) \\
                                        &  = t^\D( \varphi(f_1),\dots, \varphi(f_k))\, (\varphi^\circ(x))\\
                                      &  = t^\A( \varphi(f_1)(\varphi^\circ(x)),\dots,\varphi(f_k)(\varphi^\circ(x)) ) \\
  & = t^\A( \varphi_x(f_1(x)),\dots,\varphi_x(f_k(x))).
\end{align*}

 Finally for the bijectivity of $\varphi_x$, recall that $(\varphi^{-1})^\circ = (\varphi^\circ)^{-1}$. 
 Hence the unique mapping $(\varphi^{-1})_{\varphi^\circ(x)}$ defined by \eqref{eq:phixdef} is the inverse of $\varphi_x$.

\eqref{it:phicirc2} 
 The point $x_k$, as the limit of the net $(y_i)_{i\in I}$, is the unique cluster point of that net.
In particular, $(y_i)_{i\in I}$ has no cluster points in $X^\circ$. 
 The mapping $(\varphi^{-1})^\circ$ is continuous on $X^\circ$ by \eqref{it:phicirc1}, and therefore
 all cluster points of the net $(\varphi^\circ(y_i))_{i\in I}$ are in $\{x_1,\dots,x_n\}$.

 Next from the definition of $\varphi_x$ 
 we see that 
\[
\varphi_{y_i}(f(y_i))=\varphi(f)(\varphi^\circ(y_i))\quad
\text{for all } f\in D,\ i\in I.
\]
 Since $\lim y_i=x_k$, we have
\[
\forall f\in D\ \exists j\in I\ \forall i\geq j :\
y_i\in f^{-1}(e_k).
\]
Note that
$y_i\in f^{-1}(e_k)$ implies $\varphi_{y_i}(f(y_i))\in (\Aut\A)(e_k)$
and hence 
\[
\forall f\in D\ \exists j\in I\ \forall i\geq j:\:
\varphi(f)(\varphi^\circ(y_i)) \in (\Aut\A)(e_k).
\]
Since $\varphi$ is an automorphism of $D$, this can be written as
\[ \forall g\in D\ \exists j\in I\ \forall i\geq j:\:
g(\varphi^\circ(y_i)) \in (\Aut\A)(e_k). \]
 It follows that if $x_\ell$ is a cluster point of $(\varphi^\circ(y_i))_{i\in I}$,
 then $g(x_\ell) = e_\ell$ is in $(\Aut\A)(e_k)$.
\end{proof}

\begin{lem} \label{la:ppohom} \mbox{}
\begin{enumerate}
\item  \label{it:h0}    
 The mapping
 \[ h^\circ\colon\Aut\ABxe\rightarrow\Homeo X^\circ,\ \varphi\mapsto\varphi^\circ, \]
 is a group homomorphism.
\item \label{it:kerh0}
 For $\varphi\in\ker h^\circ$,
 \[ \varphi_*\colon X^\circ\to\Aut\A,\ x\mapsto \varphi_x, \] 
 is continuous.
\item \label{it:K0}
 The mapping
\[ p\colon\ker h^\circ \to (\Aut\A)^{X^\circ},\ \varphi\mapsto \varphi_*, \]
 is a group monomorphism.   
\end{enumerate}
\end{lem}

\begin{proof}
\eqref{it:h0} By Lemma \ref{lem:hphi}\eqref{it:phicirc1}, $\varphi^\circ\colon X^\circ\rightarrow X^\circ$ is the unique
 mapping that satisfies $\varphi(\theta_x)=\theta_{\varphi^\circ(x)}$ for all $x\in X^\circ$.
Therefore, for $\varphi_1,\varphi_2\in \Aut\ABxe$ we have
\[
\varphi_1\varphi_2(\theta_x)=\varphi_1(\theta_{\varphi_2^\circ(x)})=
\theta_{\varphi_1^\circ\varphi_2^\circ(x)},
\]
and it follows that $(\varphi_1\varphi_2)^\circ=\varphi_1^\circ\varphi_2^\circ$, as required.

 \eqref{it:kerh0}
 Let $\varphi\in\ker h^\circ$, that is, $\varphi^\circ = \id_{X^\circ}$, let $\alpha\in\Aut\A$ and $x\in X^\circ$.
 Then the definition of $\varphi_x$ in Lemma \ref{lem:hphi}\eqref{it:phix} yields
\begin{equation}\label{eq:appp}
 \alpha = \varphi_x \Leftrightarrow \alpha \pi_x = \pi_x \varphi.
\end{equation}
 Consequently
\begin{align*}
\varphi_*^{-1}(\alpha)
& = \{x \in X^\circ \st \varphi_x = \alpha\}\\
& = \{x \in X^\circ \st \alpha \pi_x = \pi_x \varphi\}\\
& = \{x \in X^\circ \st \alpha f(x) = \varphi(f)(x) \text{ for all } f\in D\} \\
&=
\Bigl[
\bigcap_{f \in D}
\bigcup_{a \in A} \bigl(
(\alpha f)^{-1}(a) \cap (\varphi(f))^{-1}(a) \bigr)
\Bigr]
\cap X^\circ .
\end{align*}
As $\alpha f$ and $\varphi(f)$ are both continuous on $X$, all the sets
$(\alpha f)^{-1}(a)$, $ (\varphi(f))^{-1}(a)$ are clopen in $X$.
It follows that $\varphi_*^{ -1}(\alpha)$ is closed in $X^\circ$. But since this is true for all $\alpha\in\Aut\A$,
the set $\varphi_*^{ -1}(\alpha)$ is in fact clopen, completing the proof that
$\varphi_*$ is continuous.

 \eqref{it:K0}
 Let $x\in X^\circ$ and $\varphi,\psi\in\ker h$ .
  Using~\eqref{eq:appp} for these maps as well as for their composition we obtain
 \[
 \varphi_x\psi_x\pi_x=\varphi_x\pi_x\psi=\pi_x\varphi\psi = (\varphi\psi)_x\pi_x.
\]
 Hence $(\varphi\psi)_x=\varphi_x\psi_x$ for all $x\in X^\circ$, that is, $p$ is a homomorphism.

 To see that $p$ is injective, let $\varphi\in\ker p$, that is, $\varphi^\circ=\id_{X^\circ}$ and
 $\varphi_x = \id_A$ for all $x\in X^\circ$.
 For $f\in D$ and $x\in X^\circ$, using \eqref{eq:appp} we have
\[ f(x) = \varphi_x(f(x)) = \pi_x\varphi\pi_x^{-1}(f(x)) = \pi_x\varphi(f) = \varphi(f) (x). \]
 Hence $\varphi=\id_D$.
\end{proof}

In general the homeomorphism $\varphi^\circ$ on $X^\circ$ from Lemma \ref{lem:hphi}
does not have a continuous extension on $X$ as shown in the following.

 \begin{exa}
Let $\A$ be a  finite non-abelian simple Mal'cev algebra with an idempotent $e$,
let $\B$ be the countable atomless Boolean algebra, and let $x_0:=00\dots$, $x_1:=11\dots$ in $2^\omega$.
We will construct an automorphism $\varphi$ of the Boolean power 
$\D:= (\A^\B)^{x_0,x_1}_{e,e}$ such that the homeomorphism $\varphi^\circ$ of $X^\circ:=2^\omega\setminus\{x_0,x_1\}$
 cannot be extended to a homeomorphism of $2^\omega$.

 To this end, consider the following clopen sets in $2^\omega$:
\begin{align*}
&b_i:=\{ \sigma\in 2^\omega \st \sigma_0=i\} \quad \text{for } i\in\{0,1\},\\
&b_{ij}:=\{\sigma\in 2^\omega \st \sigma_l=i \text{ for } l\in [j], \sigma_{j}=1-i\} \quad \text{for } i\in\{0,1\}, j\geq 1.
\end{align*}
Note that
\[
(\bigcup_{j\geq 1} b_{ij} )\cup \{x_i\} =b_i \quad \text{for } i\in\{0,1\}.
\]
 Define $\psi: X^\circ\rightarrow X^\circ$ by
\begin{multline*}
\psi\colon   (\underbrace{i,\dots,i}_{j},1-i,\sigma_{j+2},\sigma_{j+3},\dots)\\
\mapsto
\begin{cases}
(\underbrace{i,\dots,i}_{j},1-i,\sigma_{j+2},\sigma_{j+3},\dots) &\text{if } j \text{ is odd},\\
(\underbrace{1-i,\dots,1-i}_{j},i,\sigma_{j+2},\sigma_{j+3},\dots)&\text{if } j \text{ is even}.
\end{cases}
\end{multline*}
In particular,
\[
\psi (b_{ij})=\begin{cases} b_{ij} &\text{if } j \text{ is odd},\\
b_{1-i,j} &\text{if } j \text{ is even}.
\end{cases}
\]
It is clear that $\psi$ is indeed a homeomorphism of $X^\circ$.
But $\psi$ cannot be extended to a homeomorphism of $2^\omega$. Indeed, supposing it could, 
we would have
\begin{align*}
\psi(b_0) & =\psi(\bigcup_{i\in \N} b_{0i}\ \cup\: \{x_0\})
=\bigcup_{in\in\N} \psi(b_{0i})\ \cup\: \{\psi(x_0)\}
\\
&=\bigcup_{i\in \N} b_{0,2i-1}\ \cup\ \bigcup_{i\in\N} b_{1,2i}\ \cup\: \{\psi(x_0)\}.
\end{align*}
The set $b_0$ is closed in $2^\omega$, whereas the set $\psi(b_0)$ has both $x_0$ and $x_1$ in its closure,
but contains only one of them, namely $\psi(x_0)$.

Now define $\varphi\colon D\rightarrow D$ by
\begin{equation}
\label{eq:defphiex}
(\varphi(f))(x):= \begin{cases} f\psi^{-1}(x)&\text{if } x\in X^\circ,\\
e&\text{if } x\in\{x_0,x_1\}.\end{cases}
\end{equation}
We show that $\varphi$ is an automorphism of $\D$ and that ${\varphi^\circ=\psi}$.

To verify that $\varphi$ is well-defined, we show that $\varphi(f)$ is continuous on $2^\omega$.
It is certainly continuous on $X^\circ$ as a composition of two continuous maps.
We now claim that $\varphi(f)$ is also continuous at $x_i$ for $i\in\{0,1\}$.
Since $f\in D$, the set $c:=f^{-1}(e)$ is clopen in $X$ and contains $x_0,x_1$.
Therefore there exists $k\in\N$ such that
$d:=\bigcup_{j\geq k} (b_{0j}\cup b_{1j}) \subseteq c$.
 Since $\psi(d)=d$, it follows that 
\[
d = \psi^{-1}(d) \subseteq \psi^{-1}(c) = \psi^{-1}(f^{-1}(e))
= (f \psi)^{-1}(e) = (\varphi(f))^{-1}(e).
\]
Hence 
\[
(\varphi(f))^{-1}(e)= \{ x_0,x_1\} \mathbin{\dot{\cup}} d\mathbin{\dot{\cup}} \bigcup_{i=0}^1\bigcup_{j=1}^{k-1} 
\bigl(b_{ij} \cap (\varphi(f))^{-1}(e)\bigr).
\]
Note that 
\[
\{ x_0,x_1\} \cup d= X\setminus (\bigcup_{i=0}^1\bigcup_{j=1}^{k-1} 
b_{ij})
\] 
is clopen.
Also for $j\in [k-1]$ we have
\begin{align*}
b_{ij} \cap (\varphi(f))^{-1}(e)
&= \{x \in b_{ij} \mid \varphi(f)(x) = e\} \\ & = \{x \in b_{ij} \mid f\psi^{-1}(x) = e\} \\
&=
\begin{cases}
b_{ij} \cap f^{-1}(e), & j \text{ odd}, \\
\psi(b_{1-i,j} \cap f^{-1}(e)), & j \text{ even},
\end{cases}
\end{align*}
and hence each $b_{ij} \cap (\varphi(f))^{-1}(e)$ is clopen.
Therefore $(\varphi(f))^{-1}(e)$ is clopen and so $\varphi(f)\in D$ as required.
It is now straightforward to show that $\varphi$ is bijective and also that it is a homomorphism because of the
componentwise definition of operations in $\D$. Finally from \eqref{eq:defphiex} it easily follows that
$\varphi^\circ=\psi$, a homeomorphism of $X^\circ$
that cannot be extended to a homeomorphism of $2^\omega$.
\end{exa}

 In contrast to the previous example we show in the following theorem that $\varphi^\circ$ can be extended to $X$
 for automorphisms $\varphi$ of $\ABxe$ if $e_1,\dots,e_n$ are in distinct orbits under $\Aut\A$.
 Under this condition we can even give a semidirect decomposition of the automorphism group of $\ABxe$.
 
\begin{thm} \label{thm:AutABxe}
 Let $\A$ be a finite simple non-abelian Mal'cev algebra with idempotents $e_1,\dots,e_n$
 in distinct $\Aut\A$-orbits,
 let $\B$ be a Boolean algebra with Stone space $X$,  and let $x_1,\dots,x_n\in X$ be distinct.
 For $\varphi\in\Aut \ABxe$ and $x\in X$, define
 $\overline{\varphi}(x) := \begin{cases}  x & \text{if } x\in\{x_1,\dots,x_n\}, \\ 
 \varphi^\circ(x) & \text{else}. \end{cases}$ 
 \begin{enumerate}
 \item\label{it:AAB1}
The mapping
\[ h\colon \Aut \ABxe \to (\Homeo X)_{\{x_1,\dots,x_n\}},\ \varphi \mapsto \overline{\varphi}, \]
 is a group epimorphism and   
\[ \Aut \ABxe \cong \ker h \rtimes (\Homeo X)_{\{x_1,\dots,x_n\}}. \]
\item\label{it:AAB2}
  Let $K\subseteq (\Aut\A)^{X^\circ}$ be the set of all continuous maps $\psi\colon X^\circ\to \Aut\A$ such that
  $\psi^{-1}((\Aut\A)_{e_i}) \cup \{x_i\}$
 is open in $X$ for each $i\in [n]$.
 Then
\[ p\colon \ker h\to K,\ \varphi\mapsto \varphi_*, \] 
is a group isomorphism.
\end{enumerate}
\end{thm}

 We note that $K$ in Theorem~\ref{thm:AutABxe}\eqref{it:AAB2} can be viewed as the closure of
 $\left((\Aut\A)^\B\right)^{x_1,\dots,x_n}_{1,\dots,1}$ in $\Aut\ABxe$ under the topology of pointwise convergence.
 For $\A$ a group and $n=1$ this was already observed by Bryant and Evans~\cite[Section 2]{BE:SIP}.

\begin{proof}[Proof of Theorem~\ref{thm:AutABxe}]
\eqref{it:AAB1}
 To see $\overline{\varphi}\in(\Homeo X)_{\{x_1,\dots,x_n\}}$, first recall that the restriction
 $\varphi^\circ$ of $\overline{\varphi}$ is a homeomomorphism of $X^\circ$ by Lemma~\ref{lem:hphi}\eqref{it:phicirc1}.
 Let $(y_i)_{i\in I}$ be a convergent net in $X^\circ$ with limit $x_k$ for $k\in[n]$. Then
 $x_k$ is the unique cluster point (i.e. the limit) of  $(\varphi^\circ(y_i))_{i\in I}$	
 by Lemma~\ref{lem:hphi}\eqref{it:phicirc2} and the assumption that $e_1,\dots,e_n$ are in distinct $\Aut\A$-orbits.
 Hence $\overline{\varphi}$ is continuous on $X$. 
 Since $\overline{\varphi}$ fixes $x_1,\dots,x_n$ by definition and is bijective on $X^\circ$ by Lemma \ref{lem:hphi}\eqref{it:phicirc1}, it follows that $\overline{\varphi}\in(\Homeo X)_{\{x_1,\dots,x_n\}}$.

 Next note that $h$ is a group homomorphism by Lemma~\ref{la:ppohom}\eqref{it:h0}.
 That $h$ is surjective and that $\ker h$ has a complement follows from the claim that
\begin{equation} \label{eq:hg} 
hg \text{ is the identity on } (\Homeo X)_{\{x_1,\dots,x_n\}},
\end{equation}
with $g$  defined in~\eqref{eq:defg}.
 To see this let $\D := \ABxe$ and let $\psi\in (\Homeo X)_{\{x_1,\dots,x_n\}}$. Then $hg(\psi) = \overline{\psi^\D}$.
 Clearly $\overline{\psi^\D}(x_i) = x_i$ for all $i\in [n]$. Else for $x\in X^\circ$ we note that
\begin{align*}
  \theta_{\overline{\psi^\D}(x)} & = \psi^\D(\theta_x) \\
                                 & = \psi^\D( \{ (f,g)\in D^2 \st f(x) = g(x) \}) \\
   & = \{ (f\psi^{-1},g\psi^{-1})\in D^2 \st f(x) = g(x) \} \\
   & = \{ (f',g')\in D^2 \st f'\psi(x) = g'\psi(x) \} \\
   & = \theta_{\psi(x)}.
\end{align*}
 Thus $\overline{\psi^\D} = \psi$ and~\eqref{eq:hg} is proved.

 \eqref{it:AAB2}
 Note that $\ker h = \ker h^\circ$ for $h^\circ$ as in Lemma~\ref{la:ppohom}.
 Hence by that lemma it only remains to show that
\begin{equation} \label{eq:pkh}
  p(\ker h) = K.
\end{equation}  
 For $\subseteq$, let $\varphi\in\ker h$.
 We show that $\varphi_*^{-1}((\Aut\A)_{e_i}) \cup \{x_i\}$ is also open in $X$ for all $i\in [n]$.
 Consider an arbitrary clopen $b$ in $X$ such that $b\cap \{x_1,\dots,x_n \} = x_i$.
Let  $f_{b}\in D$ be an arbitrary element satisfying $f_{b}(b) = e_i$.
Using \eqref{eq:appp}, for any $x\in b\setminus \{x_i\}$ we have $\varphi_x\in(\Aut\A)_{e_i}$ if and only if
$\pi_x\varphi(f_b) = e_i$. The latter is equivalent to $x\in (\varphi(f_b))^{-1}(e_i)$.
 Hence
\[ \bigl[\varphi_*^{-1}((\Aut\A)_{e_i})\cup \{x_i\}\bigr] \cap b = (\varphi(f_b))^{-1}(e_i) \cap b\]
 is clopen in $X$. 
 Note that $b$ can be chosen to contain any given $x\in X^\circ$, and hence we have
\begin{align*}
   & \varphi_*^{-1}((\Aut\A)_{e_i}) \cup \{x_i\} =  \\
   & \bigcup\Bigl\{ \bigl[\varphi_*^{-1}((\Aut\A)_{e_i}) \cup \{x_i\}\bigr] \cap b \st b \text{ clopen in } X,
   b\cap \{x_1,\dots,x_k \} = \{x_i\} \Bigr\},
\end{align*}
 a union of open sets in $X$.  
 This completes the proof that ${p(\ker h)\subseteq K}$.
 
 For the converse inclusion, let $\psi\in K$.
 For $f\in D$ define
\[ \psi'(f)\colon X \to A,\ x\mapsto \begin{cases} [\psi(x)](f(x)) & \text{if } x\in X^\circ, \\ f(x) & \text{else}. \end{cases} \]
 Clearly $\psi'$ is the only possible choice for the preimage of $\psi$ under $p$.
 But it remains to check that $\psi'$ actually maps $D$ to itself.
 To see this let $f\in D$. From the definition we have $[\psi'(f)](x_i) = f(x_i) = e_i$ for all $i\in [n]$.  
 To show that $\psi'(f)$ is continuous, let $a\in A$ and consider
\[ X^\circ \cap [\psi'(f)]^{-1}(a) = \bigcup_{\alpha\in\Aut\A} \{ x\in X^\circ \st \psi(x) = \alpha, f(x) = \alpha^{-1}(a) \}. \]
 Since $\psi$ and $f$ are continuous on $X^\circ$ and $X$, respectively, and $X^\circ$ is open in $X$,
 we obtain that
\begin{equation} \label{eq:psi'}
 X^\circ \cap [\psi'(f)]^{-1}(a) \text{ is open in } X \text{ for all } a\in A.
\end{equation}
 In particular, if $a\not\in\{e_1,\dots,e_n\}$, then $[\psi'(f)]^{-1}(a) = X^\circ \cap [\psi'(f)]^{-1}(a)$ is open in $X$.

 On the other hand, for $i\in [n]$, the set $b := f^{-1}(e_i)$ is open in $X$ and contains $x_i$.
 By assumption $c := \psi^{-1}((\Aut\A)_{e_i}) \cup \{x_i\}$ is open in $X$ as well. Hence $b\cap c$ is open and contains
 $x_i$.
 Since $\psi'(f)(b\cap c) = \{e_i\}$, we obtain with~\eqref{eq:psi'} that
\[ [\psi'(f)]^{-1}(e_i) = [X^\circ \cap [\psi'(f)]^{-1}(e_i)] \cup [b\cap c] \]
 is open in $X$.
 Thus $\psi'(f)$ is continuous and consequently in $D$. 

 From the pointwise definition of $\psi'$ it is straightforward that $\psi'$ is an automorphism of $\D$,
 $\psi'\in\ker h$ and $p(\psi') = \psi$.
 Thus $p(\ker h) \supseteq K$ and \eqref{eq:pkh} is proved. Hence $p$ is an isomorphism as required.
\end{proof}

\subsection{Isomorphisms}

In this subsection we discuss isomorphisms between filtered Boolean powers, in the case where $\B$ is the countable atomless Boolean algebra. 
The aim is to show that in a Boolean power $\D=\ABxe$ we can add or remove filtering idempotents which
  belong to the orbits of the remaining ones 
  under $\Aut\A$ without changing the isomorphism type of $\D$.

 As a stepping stone towards the main result of this section, the following lemma establishes one important special case when
the number of filter points in a filtered Boolean power can be reduced.   

\begin{lem}  \label{lem:isoproduct}
  Let $\A$ be an algebra with idempotent $e$,
  let $\B$ be the countable atomless Boolean algebra, and let $x_0\in 2^\omega$.

  Then $(\AB)^{x_0}_{e} \times (\AB)^{x_0}_{e} \cong (\AB)^{x_0}_{e}$.
\end{lem}

\begin{proof}
Let $X_1,X_2$ be two disjoint copies of $2^\omega$ and fix $x_i\in X_i$ for $i\in [2]$.
 We will represent the two copies of $(\AB)^{x_0}_{e}$ as
\[
 D_i:= \{ f\colon X_i\rightarrow A \st f\text{ is continuous and } f(x_i)=e\} \quad \text{ for } i\in [2].
\]
Consider now the space $X_1\cup X_2$ with the union topology; it is homeomorphic to $2^\omega$.
Let $\sim$ be the equivalence relation on $X_1\cup X_2$ with the equivalence classes
$\{ x_1,x_2\}$ and $\{x\}$ for $x\in (X_1\cup X_2)\setminus \{x_1,x_2\}$.
The quotient space $X :=(X_1\cup X_2)/{\sim}$ is again homeomorphic to $2^\omega$.
This is easy to show directly by verifying that it is Hausdorff, compact, has no isolated points and has a countable basis
of clopen sets. Alternatively, one can note that $\sim$ is a Boolean equivalence in the sense of \cite{Mo:HBA1},
and that the resulting quotient $(X_1\cup X_2)/\!\sim$ has a countable basis of clopen sets and no isolated points.
 We will take its points to be $(X_1\setminus\{x_1\})\cup (X_2\setminus\{x_2\})\cup \{x_{12}\}$ and note that there are
 two types of open sets:
\begin{itemize}
\item
$U_1\cup U_2$ where $ U_i$ is open in  $X_i$ and  $x_i\not\in U_i$  for  $i\in[2]$;
\item
$(U_1\setminus\{x_1\} \cup (U_2\setminus\{ x_2\})\cup\{x_{12}\}$ where $U_i$  is open in $X_i$ and $x_i\in U_i$ for $i\in[2]$.
\end{itemize}
Let
\[
D:=\{ f \colon X\rightarrow X \st f\text{ is continuous and } f(x_{12})=e\}
\]
be yet another copy of $(\AB)^{x_0}_{e}$. We will prove the lemma by showing that
$\psi\colon D_1\times D_2\rightarrow D$  defined by
\[
\psi(f_1,f_2)(x):=\begin{cases} f_1(x) & \text{if } x\in X_1\setminus\{x_1\},\\
f_2(x) & \text{if } x\in X_2\setminus\{x_2\},\\
e & \text{if } x=x_{12},\end{cases}
\]
is an isomorphism.

We prove that $\psi$ is well-defined, i.e. that $\psi(f_1,f_2)\in D$ for all $f_i\in D_i$, $i\in [n]$.
Let $a\in A$.
If $a\neq e$, then 
\[
(\psi(f_1,f_2))^{-1}(a)=f_1^{-1}(a)\cup f_2^{-1}(a),
\]
where each $f_i^{-1}(a)$ is a clopen in $X_i$ not containing $x_i$.
Hence $(\psi(f_1,f_2))^{-1}(a)$ is clopen in $X$.
Similarly, for $a=e$ we have 
\[
(\psi(f_1,f_2))^{-1}(a)=(f_1^{-1}(e)\setminus\{x_1\})\cup (f_2^{-1}(e)\setminus \{x_2\})\cup \{x_{12}\},
\]
again a clopen in $X$.
Hence indeed $\psi(f_1,f_2)\in D$.
It is clear from the definition that $\psi$ is bijective. It is a homomorphism by the componentwise definition of
operations in $D_1$, $D_2$ and $D$.
\end{proof}  

We now show that any filtered Boolean power of any algebra by the countable atomless Boolean algebra is isomorphic to a Boolean power in which the filtering idempotents are in distinct $\Aut \A$-orbits.

\begin{cor} \label{cor:isoreduced}
 Let $\A$ be a finite algebra, $n\leq m$
 and let $e_1,\dots,e_n$ be orbit representatives of
 idempotent elements $e_1,\dots,e_m$ in $\A$ under $\Aut\A$.
 Let $\B$ be the countable atomless Boolean algebra, and let $x_1,\dots,x_m\in 2^\omega$ be distinct.

 Then $(\A^\B)^{x_1,\dots, x_m}_{e_1,\dots,e_m}$ is isomorphic to $\ABxe$.
\end{cor}

\begin{proof}
  By partitioning $2^\omega$ into disjoint clopens $b_1,\dots,b_m$ with $x_i\in b_i$,
  and recalling that all $b_i$ are homeomorphic to $2^\omega$, we obtain
\begin{equation} \label{eq:prod}
(\A^\B)^{x_1,\dots, x_m}_{e_1,\dots,e_m}\cong \prod_{i=1}^m (\A^\B)^{x_0}_{e_i}
\end{equation}
 where $x_0\in 2^\omega$.

Now, for every $j\in \{ n+1,\dots,m\}$ there exist $i\in [n]$ and $\alpha\in\Aut \A$ such that
$\alpha(e_j)=e_i$.
Clearly $(\A^\B)^{x_0}_{e_j}\rightarrow (\A^\B)^{x_0}_{e_i},\ f\mapsto \alpha f,$ is an isomorphism.
Together with Lemma \ref{lem:isoproduct} this yields 
\[
(\A^\B)^{x_0}_{e_i}\times (\A^\B)^{x_0}_{e_j}\cong (\A^\B)^{x_0}_{e_i}\times (\A^\B)^{x_0}_{e_i}\cong
(\A^\B)^{x_0}_{e_i}.
\]
The result follows from~\eqref{eq:prod}.
\end{proof}

\subsection{$\omega$-categoricity} \label{sec:omegacat}

 Macintyre and Rosenstein's characterization of $\omega$-categorical filtered Boolean powers in~\cite{MR:CRN}
 yields the following special case and in particular Theorem~\ref{thm:SIP}\eqref{it:omegacat}.

 \begin{thm}   \label{thm:omegacat}
 Let $\A$ be a finite algebra, let $n \geq 0$, let $e_1 ,\dots , e_n$ be idempotents of $\A$,
 let $\B$ be the countable atomless Boolean algebra, and let $x_1 ,\dots, x_n \in 2^\omega$ be distinct.

 Then $\ABxe$ is $\omega$-categorical. 
\end{thm}

\begin{proof}
Since Macintyre and Rosenstein use Stone spaces of maximal ideals instead of ultrafilters in~\cite{MR:CRN}, 
 we change to that dual perspective in the following.
 Note that the complement of an ultrafilter $x'_i := B\setminus x_i$ is a maximal ideal of $\B$ for $i\in [n]$.
 We can identify $x'_i$ with the set of all clopens in the Stone space of $\B$ that do not contain the point $x_i$.

Recall that the Heyting algebra of ideals of $\B$ has operations
\[ I\wedge J := I\cap J,  \quad I\vee J := I+J, \quad I\rightarrow J := \{b\in B \st bI \subseteq J \} \] 
 for ideals $I,J$. For $U,V\subseteq [n]$, it is straightforward to check that
\begin{align*}
  & \bigwedge_{i\in U} x'_i   \wedge  \bigwedge_{i\in V} x'_i  = \bigwedge_{i\in U\cup V} x'_i, \\
  & \bigwedge_{i\in U} x'_i   +  \bigwedge_{i\in V} x'_i  = \bigwedge_{i\in U\cap V} x'_i, \\
  & \bigwedge_{i\in U} x'_i \rightarrow  \bigwedge_{i\in V} x'_i  = \bigwedge_{i\in V\setminus U} x'_i.
\end{align*}
 So the Heyting algebra $\H_0$ of ideals of $\B$ generated by $x'_1,\dots,x'_n,0,B$ has universe
\[ \{ \bigwedge_{i\in U} x'_{i} \st U\subseteq [n] \} \cup \{0\}. \]
In particular $\H_0$ is finite and $\B/J$ has finitely many atoms for any ideal $J$ in $\H_0$.
 Now the augmented Boolean algebra $(\B,x'_1,\dots,x'_n)$ is a first order structure
 with operations from the Boolean algebra and unary relations that are ideals $x'_1,\dots,x'_n$ (see \cite[p. 137]{MR:CRN}).
 Then $(\B,x'_1,\dots,x'_n)$ is $\omega$-categorical by~\cite[Theorem~7]{MR:CRN}.
 Hence $\ABxe$ is $\omega$-categorical by~\cite[Theorem~5]{MR:CRN}.
\end{proof}

\section{The  generalized Fra{\"i}ss{\'e} limit} \label{sec:Flim}

We apply standard model theoretic arguments to show that the generalized Fra{\"i}ss{\'e} limit of
$\{\mathbf{A}^k \st k\geq 1\}$ exists and is isomorphic to a filtered Boolean power.

\begin{proof}[Proof of Theorem~\ref{thm:Flim}\eqref{it:Flim}]
 Let $\A$ be a finite simple non-abelian Mal'cev algebra. 
We prove that the generalized Fra{\"i}ss{\'e} limit of $K := \{\mathbf{A}^k \st {k\geq 1}\}$ exists by using 
Lemma \ref{la:genFr}. For this we need to check that all $\A^k$ are finitely generated, 
that $K$ is countable, and that it has the JEP and AP. The first three statements are clear.

 To see that $K$ has the AP, consider $u,v,w\in\N$ and embeddings $\varphi\colon \A^u\to\A^v$, $\psi\colon \A^u\to\A^w$.
 Since $\A$ is simple, we can identify the endomorphisms of $\A$ that are not automorphisms with the set of idempotents
 $E := \{e_1,\dots,e_n\}$ of $\A$. 
 From the Foster--Pixley Theorem (Lemma \ref{la:FP}) it follows that for every $i\in [v]$ there exist
 $\alpha_i\in\Aut\A \cup E$ and $j_i\in [u]$ such that
\[ \varphi(a_1,\dots,a_u) = (\alpha_1(a_{j_1}),\dots,\alpha_v(a_{j_v})) \text{ for all } a_1,\dots,a_u\in A. \]
 Further $\{j_i\st \alpha_i\in\Aut\A\} = [u]$ since $\varphi$ is an embedding.
By permuting the copies of $\A$ and renaming their elements we may assume that
\[ \varphi(a_1,\dots,a_u) = (\underbrace{a_1,\dots,a_1}_{p_1 \text{ times}},\dots,\underbrace{a_u,\dots,a_u}_{p_u \text{ times}},\underbrace{e_1,\dots,e_1}_{q_1 \text{ times}},\dots,\underbrace{e_n,\dots,e_n}_{q_n \text{ times}}), \]
 with multiplicities $p_1,\dots,p_u\geq 1, q_1,\dots,q_n\geq 0$.
Analogously,
\[ \psi(a_1,\dots,a_u) = (\underbrace{a_1,\dots,a_1}_{r_1 \text{ times}},\dots,\underbrace{a_u,\dots,a_u}_{r_u \text{ times}},\underbrace{e_1,\dots,e_1}_{s_1 \text{ times}},\dots,\underbrace{e_n,\dots,e_n}_{s_n \text{ times}}), \]
 where $r_1,\dots,r_u\geq 1, s_1,\dots,s_n\geq 0$. 
Letting
 $v_i:=\max(p_i,r_i)$ for $i\in [u]$, $w_i:= \max(q_i,s_i)$ for $i\in [n]$ and
 $m := \sum_{i=1}^u v_i+ \sum_{i=1}^n w_i$ it is now straightforward to define embeddings $\varphi'\colon\A^v\to\A^m$
 and $\psi'\colon\A^w\to\A^m$ as 
 \[ \varphi'(a_1,\dots,a_v) := (\underbrace{a_1,\dots,a_1}_{v_1-p_1+1 \text{ times}},a_2,\dots,a_{p_1},\underbrace{a_{p_1+1},\dots,a_{p_1+1}}_{v_2-p_2+1 \text{ times}},a_{p_1+2},\dots,a_{p_1+p_2},\dots), \]
 \[ \psi'(a_1,\dots,a_w) := (\underbrace{a_1,\dots,a_1}_{v_1-r_1+1 \text{ times}},a_2,\dots,a_{r_1},\underbrace{a_{r_1+1},\dots,a_{r_1+1}}_{v_2-r_2+1 \text{ times}},a_{r_1+2},\dots,a_{r_1+r_2},\dots), \]
 taking  the sequences of multiplicities of variables to be the concatenations of:
 \begin{align*}
 &(v_i-p_i+1,\underbrace{1,\dots,1}_{p_i-1})\ (i\in [u])\text{ and }
 (w_i-q_i+1,\underbrace{1,\dots,1}_{q_i-1})\ (i\in [n]) \text{ for } \varphi',\\
 &(v_i-r_i+1,\underbrace{1,\dots,1}_{r_i-1})\ (i\in [u])\text{ and }
 (w_i-s_i+1,\underbrace{1,\dots,1}_{s_i-1})\ (i\in [n]) \text{ for } \psi'.
 \end{align*}
 Then $\varphi'\varphi=\psi'\psi$ since both map $(a_1,\dots,a_u)$ to
\[  (\underbrace{a_1,\dots,a_1}_{v_1 \text{ times}},\dots,\underbrace{a_u,\dots,a_u}_{v_u \text{ times}},\underbrace{e_1,\dots,e_1}_{w_1 \text{ times}},\dots,\underbrace{e_n,\dots,e_n}_{w_n \text{ times}}). \]

 Since $K$ has the JEP and AP, by Lemma \ref{la:genFr} there exists a unique (up to isomorphism) countable algebra
 $\C$ such that
\begin{enumerate}
\item[(i)] 
 the set of finitely generated subalgebras of $\C$ is the class of all finitely generated algebras that can be embedded into
 algebras in $K$ up to isomorphism,
\item[(ii)] $\C$ is a direct limit of algebras in $K$, and
\item[(iii)] every isomorphism between subalgebras of $\C$ that are isomorphic to some element in $K$
  extends to an automorphism of $\C$ (i.e., $\C$ is $K$-\emph{homogeneous}).
\end{enumerate}

Let $\B$ be the countable atomless Boolean algebra, and let $x_1,\dots,x_n\in 2^\omega$ be distinct.
 To see that the filtered Boolean power $\D:=\ABxe$
 is isomorphic to the generalized Fra{\"i}ss{\'e} limit $\C$
 it suffices to show that $\D$ satisfies (i), (ii), (iii). 
 Properties (i) and (ii) follow from the fact that $\D$ is locally finite and that every finite subalgebra of $\D$
 is contained in a subalgebra isomorphic to some finite direct power of $\A$, as in the proof of Lemma \ref{lem:con}.
To prove (iii) we emulate \cite[Lemma 7.1.4]{Ho:MT}, and show that
 $\D$ is \emph{weakly $K$-homogeneous}, in the sense that for any $u\leq v$ and
 any embeddings $\varphi\colon \A^u \to \A^v$ and $\psi\colon \A^u \to \D$ there exists an embedding
 $\psi'\colon \A^v \to \D$ such that $\psi'\varphi = \psi$.
 That weak $K$-homogeneity implies $K$-homogeneity is straightforward using a standard back-and-forth argument as in~\cite[cf. Lemma 7.1.4]{Ho:MT}.

To verify weak $K$-homogeneity, as above we may assume that 
\[ \varphi(a_1,\dots,a_u) = (\underbrace{a_1,\dots,a_1}_{p_1 \text{ times}},\dots,\underbrace{a_u,\dots,a_u}_{p_u \text{ times}},\underbrace{e_1,\dots,e_1}_{q_1 \text{ times}},\dots,\underbrace{e_n,\dots,e_n}_{q_n \text{ times}}) \]
with multiplicities $p_1,\dots,p_u\geq 1, q_1,\dots,q_n\geq 0$ for all $a_1,\dots,a_u\in A$.
 For $\psi\colon \A^u \to \D$ there exist $m_i\in\N$ ($i\in [u]$), automorphisms 
 $\alpha_{ij}$ of $\A$ ($i\in [u], j\in [m_i]$), and  a partition of $2^\omega$ into clopens 
 $c_{ij}$ ($i\in [u], {j\in [m_i]}$) and $b_i$ ($i\in [n]$)  such that
\[ \psi(a_1,\dots,a_u)(x)= \begin{cases} \alpha_{ij}(a_i) & \text{if } x\in c_{ij} \text{ for } i\in [u], j\in [m_i], \\
    e_i & \text{if } x\in b_i \text{ for } i\in [n]. \end{cases} \]
 Clearly we may assume that $m_i \geq p_i$ for all $i\in [u]$. Now define an embedding  $\psi'\colon\A^v\to\D$ as follows.
 Subdivide each $b_i$ into $q_i+1$ clopens $b_{ij}$ ($0\leq j\leq q_i$) with $x_i\in b_{i0}$.
 Set
\begin{align*}
 p_{ij} &:=\begin{cases} \sum_{k=1}^{i-1} p_k+j &\text{if } i\in [u],\ j\in [p_i],\\
\sum_{k=1}^{i} p_k&\text{if } i\in [u],\ j\in [m_i]\setminus [p_i],\\
 \end{cases}
 \\
\intertext{and}
q_{ij}&:= \textstyle\sum_{k=1}^{u} p_k+\sum_{k=1}^{i-1} q_i+j\ \text{ if } i\in [n], j\in [q_i].
 \end{align*}
 Then define
 \[ \psi'(a_1,\dots,a_v)(x)= \begin{cases} \alpha_{ij}(a_{p_{ij}}) & \text{if } x\in c_{ij} \text{ for } i\in [u], j\in [m_i], \\
   
    a_{q_{ij}} & \text{if } x\in b_{ij} \text{ for } i\in [n], j\in [q_i], \\
    e_i & \text{if } x\in b_{i0} \text{ for } i\in [n]. \end{cases} \]
It is straightforward to verify that $\psi'$ is indeed an embedding and that $\psi'\varphi = \psi$.

This completes the proof of weak homogeneity of $\ABxe$, and hence of the theorem.
\end{proof}

 We conclude this section with examples showing that in general 
\begin{enumerate}
\item the class $K$ of finite direct powers of a finite algebra $\A$ may not have the AP and 
\item the $K$-homogeneous limit $\D$ above may fail to be homogeneous.
\end{enumerate}  

\begin{exa}
 It is well-known and not hard to show that the class of finite abelian groups has the AP. In particular, the class of
 finite direct powers of any finite abelian group has the AP.

 We claim that this does not generalize to $2$-nilpotent groups. Let $\A = D_8\times C_2$ be the direct product
 of the dihedral group $D_8 := \langle a,b \st a^4=b^2=1, a^b = a^{-1}\rangle$ of order $8$ and the cyclic group
 $C_2:=\langle -1\rangle$ of order $2$.
 We define embeddings $\varphi_1,\varphi_2\colon \A\to\A^2$ by
\[
  \varphi_1(x,(-1)^z) := \bigl((x,1),(a^{2z},1)\bigr), \quad \varphi_2(x,(-1)^z) := \bigl((x,1),(b^z,1)\bigr),
\]
 for $x\in D_8, z\in\Z$.   
 Assume that there exists a group $\B$ and embeddings $\psi_i\colon \A^2\to\B$ such that
 $\psi_1\varphi_1=\psi_2\varphi_2$.
 Note that $a^2 = [a,b]$ is a commutator and hence $\psi_1\varphi_1(1,-1)$ is an element of the derived subgroup of $\B$.
 On the other hand, $\psi_2\varphi_2(1,-1)$ does not commute with $\psi_2 \bigl((1,1),(a,1)\bigr)$ and hence is not in
 the center of $\B$. Thus $\B$ is not $2$-nilpotent, in particular, not a direct power of $\A$.
 So $\{\A^k \st k\geq 1 \}$ does not have the AP.
\end{exa}

\begin{exa}
 For a generalized Fra{\"i}ss{\'e} limit that is not homogeneous for arbitrary finite substructures,
 let $\A$ be the alternating group of degree $6$.
 Let $2^\omega=b_1\cup b_2$ be a partition of $2^\omega$ into two non-empty clopens, and let $x_1\in b_1$.
 Consider the  two elements $f,g$ of the filtered Boolean power $\D=(\A^\B)^{x_1}_1$ defined as follows:
\[
f(x):=\begin{cases} 1 & \text{if } x\in b_1\\ (1\: 2\: 3) & \text{if } x\in b_2\end{cases},\quad
g(x):=\begin{cases} 1 & \text{if } x\in b_1\\ (1\: 2\: 3)(4\: 5\: 6) & \text{if } x\in b_2\end{cases}.
\]
Each of these two elements generates a cyclic subgroup of $\D$ of order $3$.
Since $\Aut(\A)$ is the symmetric group $S_6$ acting by conjugation,  
there is no $\alpha\in \Aut(\A)$ with $\alpha(1\: 2\: 3)=(1\: 2\: 3)(4\: 5\: 6)$.
By our description of automorphisms from Theorem \ref{thm:AutABxe} it follows that
there is no automorphism of $\D$ which maps the subgroup generated by $f$ to that generated by $g$. 
\end{exa}

\section{The variety generated by a simple algebra} \label{sec:variety}

We investigate the countable free algebra in the variety $V$ generated by a finite simple non-abelian
 Mal'cev algebra $\A$.

 \begin{proof}[Proof of Theorem~\ref{thm:Flim}\eqref{it:free}]
 First we build explicit representations of finite and countable free algebras in $V$.
 For  $i\in \N$ let $x_i\colon A^\N\to A,\ a\mapsto a(i)$. 
 Then the universe of the subalgebra $\F := \langle x_i \st i\in \N \rangle$ of  $\A^{A^\N}$ generated by the $x_i$
  forms the clone of term operations on $\A$
 (see \cite[Section 4.1]{MMT:ALV1}). Hence $\F$ is free in $V$ over $\{x_i\st i\in \N\}$ (see \cite[Section 4.11]{MMT:ALV1}).
 Note that $A^\N$ with the product topology is homeomorphic to the Cantor space.
 Since $x_i$ for $i\in\N$ are continuous functions from $A^\N$ to $A$, so are all the elements of $F$.
 Hence $F$ is contained in the Boolean power $\A^\B$ for $\B$ the countable atomless Boolean algebra.

Consider the lexicographic ordering on $A^\N$, where $A$ is linearly ordered in an arbitrary way.
Let $R$ be the transversal of $\Aut\A$ orbits on $A^\N$ consisting of lexicographically minimal elements.
It is easy to see that:
\begin{itemize}
\item
For every $k\in\N$ the restriction $R|_{[k]}$
 is a set of orbit representatives of $A^k$ under $\Aut\A$.
\item
The restriction $\F|_R$ is isomorphic to $\F$.
\end{itemize}
 We will also need the following explicit description of the free algebra 
 $\F_k=\langle x_1,\dots,x_k\rangle$ of rank $k$ as a subalgebra of $\A^\B$, as well as the descriptions of the corresponding free algebras in the variety $W$ generated by the proper subalgebras of~$\A$. Let
\[ S_k := \{ a\in R \st \langle a(1),\dots,a(k)\rangle \neq A \}. \]
 Then $\F_k|_{R\setminus S_k}$ is a finite subdirect power of $\A$, hence isomorphic to a direct power of $\A$
 by the Foster--Pixley Theorem 
 (Lemma \ref{la:FP}). Clearly $\F_k$ is a subdirect product of
 its projections $\F_k|_{R\setminus S_k}$ and $\F_k|_{S_k}$. Since every quotient of the former is isomorphic to a direct power
 of $\A$ again (cf. Lemma~\ref{lem:con}) and the latter is in $W$,
 the factors $\F_k|_{R\setminus S_k}$ and $\F_k|_{S_k}$ have no non-trivial common homomorphic image.
 So Fleischer's Lemma (Lemma \ref{la:Fleischer})
 yields that $\F_k$ is the direct product of its projections to $R\setminus S_k$ and $S_k$ respectively.
 For $p \in R|_{[k]}$ let $c_p := \{ r\in R \st r|_{[k]} = p \}$ be a basic clopen in $R$, and
 note that every element of $\F_k$ is constant on every $c_p$.
 This in turn implies
\begin{equation} \label{eq:Fk}
  F_k = \left\{ f\in A^R \st \begin{array}{l}
                           f|_{S_k}\in F_k|_{S_k},\, f \text{ is constant on } c_p \\ \text{for all } p \in (R\setminus S_k)|_{[k]}
                           \end{array} \right\}.
\end{equation}
Further $\F_k|_{S_k}$ is the greatest homomorphic image of $\F_k$ that is in $W$.
  Thus $\F_k|_{S_k}$ is free in $W$ over the restrictions of $x_1,\dots,x_k$ to $S_k$.                  
Let
\[ S := \bigcap_{k\in\N} S_k = \{ a\in R \st F|_a \neq A \}. \] 
 Since $\{ a|_{[k]} \st a\in S_k \} = \{ a|_{[k]} \st a\in S \}$ for every $k\in\N$,
  the restriction of $\F_k|_{S_k}$ to $S$ is an isomorphism.
 Hence $\F_k|_S$ is free in $W$ over $\{ x_i|_S \st i\in [k]\}$ for each $k\in\N$
 and consequently $\F|_S$ is free in $W$ over $\{ x_i|_S \st i\in \N\}$.
 Since the kernel $\theta$ of the restriction map of $\F$ to $S$ is union of the kernels of the restrictions of $\F_k$
 to $S_k$ for all $k\in\N$, it follows that $\theta$ is the smallest congruence such that $\F/\theta\in W$.

 So let $e\in F$ such that $e|_{S}$ forms a singleton subalgebra of $\F|_{S}$.
 Then $e(S)= \{e_1,\dots,e_n\}$ is a set of idempotents of $\A$ containing at least one idempotent from each $\Aut\A$-orbit.
 The sets $s_j := e^{-1}(e_j)\cap S$ for $j\in [n]$ partition $S$. Also let
\[ B := \{ b \text{ clopen in } R \st b\cap s_j \in \{\emptyset, s_j\} \text{ for all } j\in [n] \}, \]
which is the universe of a Boolean algebra $\B$.

 We claim that $\B$ is the countable atomless Boolean algebra.
To see this consider an arbitrary basic clopen $c_p\in B$ where $p\in R|_{[k]}$.
Let $l>k$ and $q\in R|_{[l]}$ be such that $p$ is a subsequence of $q$ and $\langle q\rangle=\A$.
Then $c_q\cap S=\emptyset$, and hence $c_q\in B$. Also $c_q\subsetneq c_p$, and so $c_p$ is not an atom.
Thus, $\B$ is atomless, and since it is clearly countable, the claim follows.
 
 We may view $\B$ as a subalgebra of the Boolean algebra of subsets of $Y := (R\setminus S) \cup \{s_1,\dots,s_n \}$. Hence
\[ D :=  \{ f\in A^R \st f^{-1}(a)\in B \text{ for } a\in A,\ f(s_j) = \{ e_j\} \text{ for } j\in [n]  \} \]
 forms a subalgebra of $\A^R$ which is naturally isomorphic to $(\A^\B)^{s'_1,\dots,s'_n}_{e_1,\dots,e_n}$ by
 Lemma~\ref{lem:fBp}. 

 We claim that 
\begin{equation} \label{eq:etD}
 e/\theta = D.
\end{equation}
 The inclusion $\subseteq$ is immediate from the definitions and the fact that every element in $F|_R$ is continuous.
 
 For the converse, let $f\in D$. 
Let $k\in\N$ be large enough so that the following hold:
\begin{itemize} 
\item \label{it:fc} $f$ is constant on the basic clopens $c_p$ for all $p\in R|_{[k]}$,
\item \label{it:eFK} $e\in F_k$,
\item $f|_{S_k} = e|_{S_k}$.
\end{itemize}
 Then $f\in F_k$ by~\eqref{eq:Fk}. Hence~\eqref{eq:etD} is proved.

 Summing up, any subalgebra $e/\theta$ of $\F$ is isomorphic to $(\A^\B)^{s'_1,\dots,s'_n}_{e_1,\dots,e_n}$ by~\eqref{eq:etD}
 and Lemma~\ref{lem:fBp}. Here $\B$ is a countable atomless Boolean algebra, and
 the set of idempotents $e_1,\dots,e_n$ contains at least one representative from each $\Aut\A$-orbit.
 By Corollary~\ref{cor:isoreduced} we may assume that $e_1,\dots,e_n$ is the set of all idempotents of $\A$.
 Thus Theorem~\ref{thm:Flim}\eqref{it:free} is proved. 
\end{proof}

 For loops (in particular for groups) we can sharpen Theorem~\ref{thm:Flim}\eqref{it:free} as follows.
 Borrowing the notion from groups, we say a loop $\A$ is a \emph{split extension} of a normal subloop
  $N$ by a subloop $H$ if $H$ is a transversal through the cosets of $N$ in $\A$.
  We note that an analogous result holds for rings by the same proof with the obvious notational changes.

\begin{cor} \label{cor:group}
 Let $\A$ be a finite simple non-abelian loop with identity $e$,
 let $V$ be the variety generated by $\A$, and
 let $W$ be the variety generated by all proper subloops of $\A$.
\begin{enumerate}
\item \label{it:countable}
 The countable free loop $\F$ in $V$ has a normal subloop $N$ and a subloop $H$ such that:
\begin{enumerate} 
\item \label{it:kernel}
  $N\cong (\A^\B)^{x}_e$ where $\B$ is the countable atomless Boolean algebra and $x\in 2^\omega$,
\item \label{it:complement}
 $H$ is isomorphic to the countable free loop in $W$,
\item \label{it:semidirect}
$N\cap H = \{e\}$ and $F=NH$. 
\end{enumerate}
\item \label{it:V}
 Every countable loop $\C$ in $V$ is a split extension of a filtered Boolean power $(\A^\B)^x_e$
 for some Boolean algebra $\B$ by some loop in $W$.
\end{enumerate}
\end{cor}

\begin{proof}
 For~\eqref{it:countable} we reuse the explicit representation $\F = \langle x_i \st i \in\N\rangle \leq \A^R$
 from the proof of Theorem~\ref{thm:Flim}\eqref{it:free} above.

 \eqref{it:kernel}
 Again let $S := \{ a\in R \st F|_a \neq A \}$. Then
\[ N := \{ f\in F \st f(S) = e \} \]
is the kernel (understood as normal subloop, not as congruence)
of the projection $\pi_S$ from $\F$ to $\F|_S$, with the latter isomorphic to the free loop
in $W$ over
 $\{\pi_S(x_i) \st i\in\N\}$.
 Further $N$ is isomorphic to $(\A^\B)^x_e$ by Theorem~\ref{thm:Flim}\eqref{it:free}.

 \eqref{it:complement}
 Let $k\in\N$ and $S_k :=  \{ a\in R \st \langle a(1),\dots,a(k)\rangle \neq A \}$.
 Consider $\F_k := \langle x_1,\dots, x_k\rangle$.
 As in the proof of Theorem~\ref{thm:Flim}\eqref{it:free} above, we have $F_k = K_k \times H_k$, where 
  $K_k := F_k|_{R\setminus S_k}$ and  $H_k := F_k|_{S_k}$.
  Further $K_k$ is isomorphic to a finite direct power of $\A$, while $H_k$ is isomorphic to the free loop
  of rank $k$ in $W$.
The constant function $\overline{e}_k\colon R\setminus S_k \to \{e\}$ forms the unique trivial subloop
of $K_k$. So $\{\overline{e}_k\} \times H_k$ is a subloop of $\F_k$, which is isomorphic to $H_k$.
 In particular
\[ y_k\colon R\to A,\ a\mapsto \begin{cases} a(k)
  & \text{if } a\in S_k \\ e & \text{else}, \end{cases} \]
 is in $F_k$, hence in $F$.
 Then
 \[ H := \langle y_i \st i\in\N \rangle \leq \F \]
 and $H\in W$ since  $H|_a\in W$ for all $a\in R$.
 Moreover, since $\bigcap_{i\in\N} S_i = S$, it follows that $H$ is isomorphic to $\F|_S$, which in turn is isomorphic to
 the free loop in $W$ over $\{ y_i \st i\in\N \}$.

 \eqref{it:semidirect}
 Since $y_k\in H$ and $x_ky_k^{-1}\in N$, we have $x_k \in NH$ for all $k\in\N$. Thus $F = NH$.

 For proving that $N\cap H$ is trivial, let $t(y_1,\dots,y_k)\in H$ for some term $t$ and $k\in\N$ such that
 $t(y_1,\dots,y_k)|_S$ is constant $e$.
 By the latter every proper subalgebra of $\A$ satisfies the identity $t\approx e$. Since $\H$ is in $W$,
 it follows that $t(y_1,\dots,y_k) = \overline{e}$, the constant function with value $e$ on $R$. Thus $N\cap H = \{\overline{e}\}$.

 For~\eqref{it:V}, let $\C\in V$ be isomorphic to $\F/K$ for the countable free loop $\F$ in $V$ with some normal $K$.
 Let $N,H$ be as in~\eqref{it:countable}. Then $\F/K$ is a split extension of $NK/K$ by $HK/K$.
 The former is isomorphic to $(\A^{\B/\theta})^x_e$ for a congruence $\theta$ of $\B$ by Lemma~\ref{lem:con}.
 The latter is clearly in $W$.
\end{proof}

\section{Small index property} \label{sec:SIP}

\subsection{SIP for the punctured Cantor space}

 The proof of Theorem~\ref{thm:SIP}\eqref{it:SIP} requires the following generalization of Truss' result that the
 group of homeomorphisms of the Cantor space $2^\omega$ has the SIP \cite[Theorem 3.7]{Tr:IPG2}.

\begin{thm} \label{thm:Truss}
 Let $n\in\N\cup\{0\}$, let $x_1,\dots,x_n\in 2^\omega$ be distinct, let $G := (\Homeo 2^\omega)_{\{x_1,\dots,x_n\}}$,
 and let $H\leq G$ with $|G:H| < 2^{\aleph_0}$.  
 
 Then there exist $m\geq n$ and a partition of $2^\omega$ into clopen sets
 $b_1,\dots,b_m$ such that $x_i\in b_i$ for all $i\in [n]$ and $G_{\{b_1,\dots,b_m\}} \leq H$.
\end{thm}

 The proof makes up the rest of this subsection.
 Throughout,  we use the notation and assumptions of Theorem~\ref{thm:Truss}, assume $n\geq 1$ and let
$X^\circ := 2^\omega\setminus\{x_1,\dots,x_n\}$.

 The proof heavily relies on the corresponding results from \cite{Tr:IPG2} concerning the Cantor space itself (Theorem 3.7)
 and $\mathbb{Q}$ with a point removed (Theorem 2.7).
 We do not reproduce the arguments that are essentially repetitions of those from \cite{Tr:IPG2} but
 just point out the necessary adaptations.
 Reading the original proof from \cite{Tr:IPG2} verbatim, with $2^\omega$ or $\mathbb{Q}$ replaced by
 $X^\circ$, and keeping in these additional adaptations, constitutes a full proof of each such result.
 The arguments that are new to our situation are presented in full.

 Following Truss~\cite{Tr:IPG1,Tr:IPG2} we say that a subset $b$ of $2^\omega$ \emph{abuts} an element $x$ of $2^\omega$ if
 $x\not\in b$ and $b\cup\{x\}$ is closed in $2^\omega$.
 We will also talk about a sequence $(b_i)_{i\in\N}$ of subsets of $2^\omega$ \emph{converging} to a point $y\in 2^\omega$; this means that for every open neigbourhood of $y$ all but finitely many members of the sequence are fully contained in the neighborhood.
 
\begin{lem} \label{la:Xo3Tk3}\label{la:Xo4Tk3}
 Let $(b_i)_{i\in\N}$ be a sequence of pairwise disjoint clopen subsets of $X^\circ$ that converges to $x_1$
 such that either (a) each $b_i$ is clopen in $2^\omega$ or (b) each $b_i$ abuts $x_1$.
 
 Then $G_{2^\omega\setminus b_i} \leq H$ for some $i\in\N$.  
\end{lem}

\begin{proof}
 (a) This is identical to the argument in the first paragraph of the proof of~\cite[Theorem 3.7]{Tr:IPG2}.
 Here we can use that $G_{2^\omega\setminus b_i}\cong\Homeo 2^\omega$ for all $i\in\N$, which is a simple group
 by a result of Anderson~\cite{An:ASC}.

 (b) This is shown by following the proof of \cite[Lemma 2.5]{Tr:IPG2}.
 This time $G_{2^\omega\setminus b_i} \cong (\Homeo 2^\omega)_{x_1}$ for all $i\in\N$, which is not simple.
 Instead, by \cite[Theorem 3.13]{Tr:IPG1}, the latter has precisely one non-trivial proper normal subgroup
\[
 N:=\{ g\in \Homeo 2^\omega \st g \text{ fixes pointwise a neighborhood of } x_1 \}.
\]
 However, this normal subgroup has uncountable index in $(\Homeo 2^\omega)_{x_1}$, which can be proved by following
 \cite[p.\ 499, par.\ 3]{Tr:IPG2}. So as in~\cite[Lemma 2.5]{Tr:IPG2} we can show that
 $H\cap G_{2^\omega\setminus b_i} = G_{2^\omega\setminus b_i}$ for some $i\in\N$. 
\end{proof}

\begin{lem}
\label{la:Xo5Tk3}
 There exists a clopen set $b$ of $2^\omega$ such that the following hold:
 \begin{itemize}
 \item $b\cap \{x_1,\dots,x_n\} = \{x_1\}$; and
 \item $G_{2^\omega\setminus c}\leq H$ for every clopen set $c\subseteq b$ of $ 2^\omega$ not containing $x_1$.
 \end{itemize}
\end{lem}

\begin{proof}
 Suppose the statement is not true. Then we can inductively define two sequences of clopen sets $(b_i)_{i\in\N}$ and
 $(c_i)_{i\in\N}$ of $2^\omega$ as follows.
 For $i\in\N$ let $b_i$ be a clopen satisfying the following properties:
\begin{itemize}
\item
$b_i \cap \{x_1,\dots,x_n\} = \{x_1\}$,
\item
$b_i$ is contained in the ball of radius $1/i$ around $x_1$ with respect to any metric on $2^\omega$ which induces the same topology,
\item
$b_i$ is disjoint from each of $c_1,\dots ,c_{i-1}$.
\end{itemize}
By our assumption we have a clopen $c_i\subseteq b_i$ not containing $x_1$ such that $G_{2^\omega\setminus c_i}\nleq H$.
The sequence $(c_i)_{i\in\N}$ consists of pairwise disjoint clopen sets not containing $x_1$, which converges to $x_1$,
and such that
$G_{2^\omega\setminus c_i}\nleq H$ for all $i\in\N$. This contradicts Lemma \ref{la:Xo3Tk3} (a).
\end{proof}

\begin{lem}
\label{la:Xo6Tk3}
 There exists a clopen set $b$ of  $2^\omega$ such that the following hold:
\begin{itemize}
\item
$b\cap \{x_1,\dots,x_n\} = \{x_1\}$; and
\item
$G_{2^\omega\setminus b}\leq H$.
\end{itemize}
\end{lem}

\begin{proof}
 Let $b$ be the clopen set of $2^\omega$ guaranteed by Lemma \ref{la:Xo5Tk3}: thus $b\cap \{x_1,\dots,x_n\} = \{x_1\}$ and  
 for every clopen $c\subseteq b\setminus\{x_1\} =: b^\circ$ of $2^\omega$ we have $G_{2^\omega\setminus c}\leq H$.
 Let
\[
\Gamma:= \{ d \st d\subseteq b \text{ is clopen in } X^\circ \text{ abutting } x_1 \text{ and } G|_{2^\omega\setminus d} \leq H\}.
\]
 We will show that $b^\circ\in \Gamma$ via a sequence of claims.

\begin{claim}
\label{cl:Xo62Tk3}
If $c$ and $d$ are clopen sets of $X^\circ$ such that $c$, $d$ and $c\cap d$ abut $x_1$, then
$G_{2^\omega\setminus (c\cup d)} =\langle G_{2^\omega\setminus c},G_{2^\omega\setminus d}\rangle$.
\end{claim}

\begin{proofclaim}
This is identical to \cite[Lemma 2.6]{Tr:IPG2}.
\end{proofclaim}

\begin{claim}
\label{cl:Xo63Tk3}
If $c,d\in\Gamma$, then $c\cup d\in\Gamma$.
\end{claim}

\begin{proofclaim}
 This can be proved by following the paragraph preceding \cite[Lemma 2.9]{Tr:IPG2} using
 Lemma \ref{la:Xo4Tk3} (b) and Claim \ref{cl:Xo62Tk3}.
\end{proofclaim}

\begin{claim}
\label{cl:Xo64Tk3}
 If $d\in\Gamma$ and $c\subseteq b^\circ$ clopen in $2^\omega$, then $c\cup d\in\Gamma$.
\end{claim}

\begin{proofclaim}
This follows \cite[Lemma 2.9]{Tr:IPG2}.
\end{proofclaim}

The proof of Lemma \ref{la:Xo5Tk3} is now completed by following
 the Conclusion of the proof of  Theorem 2.7 in \cite[p. 501, 502]{Tr:IPG2}.
\end{proof}

\begin{proof}[Proof of Theorem~\ref{thm:Truss}]
 We use induction on $n$.
 The base case for $n=0$ is Truss' result~\cite[Theorem 3.7]{Tr:IPG2}.
 So assume $n\geq 1$ in the following.
 By Lemma \ref{la:Xo6Tk3} there exists a clopen $b_1$ in $2^\omega$ such that $b_1\cap\{x_1,\dots,x_n\} = \{x_1\}$ and
 $G_{2^\omega\setminus b_1}\leq H$.
 If $b_1= 2^\omega$, then $G=H$ and there is nothing left to prove. So suppose $b_1\neq 2^\omega$.
 Since $H\cap G_{b_1}$ has index $<2^{\aleph_0}$ in $G_{b_1} \cong (\Homeo 2^\omega_{x_2,\dots,x_n}$,
 we have a partition of $2^\omega\setminus b_1$ into clopen sets $b_2,\dots,b_m$ such that
 $(G_{b_1})_{\{b_2,\dots,b_m\}} \leq H\cap G_{b_1}$ by the induction assumption.  
 But then $H$ contains the direct product
\[ G_{2^\omega\setminus b_1}\ (G_{b_1})_{\{b_2,\dots,b_m\}} = G_{\{b_1,b_2,\dots,b_m\}} \]
 as required.
\end{proof}

\subsection{Orbits of clopens in the punctured Cantor space}
As above, let 
$x_1,\dots,x_n\in 2^\omega$ and
$X^\circ := 2^\omega\setminus\{x_1,\dots,x_n\}$.
 For any set $b\subseteq  2^\omega$ we write $\cl(b)$ for its \emph{closure} in  $2^\omega$.
 For a clopen $b$ in $X^\circ$ let
\begin{align*}
  & I := \{i\in [n] \st x_i\in \cl(b) \}, \\
  & I' := \{ i\in [n] \st x_i\in \cl(X^\circ\setminus b)  \}.
\end{align*} 
 We call the pair  $(I,I')$ the \emph{type} of $b$. Note that $I\cup I' = [n]$ but that union is not necessarily disjoint.
 
Let us partition the Cantor space into the disjoint union $2^\omega=c_1\cup\dots\cup c_n$ of disjoint clopens with $x_i\in c_i$ for all $i$. For any $i\in [n]$ and any $j\in \N$ let $c_{ij}$ be the intersection of $c_i$ and the ball of radius $1/j$ around $x_i$.
  The sets $c_{ij}$ are clopen, and satisfy $c_i=c_{i1}\supseteq c_{i2}\supseteq \dots \supseteq \{x_i\}$,
  and $\lim_{j\rightarrow\infty}c_{ij}=x_i$ for all $i$.
The sets $c_{i,j+1}\setminus c_{ij}$ for $i\in [n]$ and $j\in\N$ are pairwise disjoint, clopen, contained in $X^\circ$, and satisfy 
$\lim_{j\rightarrow\infty}c_{i,j+1}\setminus c_{ij}=x_i$ for all $i$.

 Now consider again an arbitrary clopen $b$ in $X^\circ$.
 Intersecting $b$ with each  $c_{i,j+1}\setminus c_{ij}$ and disregarding all the empty sets, yields a partition of $b$ into the disjoint union of non-empty clopens $b_{ij}$ of $2^\omega$ with $i\in [n]$ and $j\in J_i$, where $J_i$ is either empty, or a finite interval $[m_i]$ or equal to $\N$.
The case where $J_i$ is finite arises when $b$ is disjoint from 
$c_{i,j+1}\setminus c_{ij}$ for sufficiently large $j$, i.e. when $x_i\not\in\cl(b)$.
It follows that all $J_i$ are finite precisely when $b$ is closed in $2^\omega$, i.e. it is a clopen of $2^\omega$.
We can perform a similar analysis on the set $X^\circ\setminus b$,
and thus conclude that the following hold:

\begin{itemize}
\item
$I=\emptyset$ if and only if $b$ is clopen in $2^\omega$;
\item
$I'=\emptyset$ if and only if $X^\circ\setminus b$ is clopen in $2^\omega$;
\item
$I,I'\neq\emptyset$ if and only if there exist pairwise disjoint clopens $b_{ij}$ ($i\in I$, $j\in\N$) and $b_{ij}'$ ($i\in I'$, $j\in\N$) such that
\begin{equation}
\label{eq:typepart}
\setlength{\arraycolsep}{2pt}
\begin{array}{rlrl}
b&=\displaystyle\bigcup_{i\in I,\,j\in\N} b_{ij},\quad & \displaystyle\lim_{j\to \infty} b_{ij}&=x_i\ (i\in I),\\
X^\circ\setminus b&=\displaystyle\bigcup_{i\in I',\,j\in\N} b_{ij}',&\displaystyle\lim_{j\to \infty} b_{ij}'&=x_i\ (i\in I').
\end{array}
\end{equation}
\end{itemize}

 We can characterize the orbits of clopens in $X^\circ$ under $(\Homeo 2^\omega)_{\{x_1,\dots,x_n\}}$ as follows.

\begin{lem} \label{lem:types}
 For $I,I'\subseteq [n]$ with $I\cup I'=[n]$, the set of clopens in $X^\circ$ of type $(I,I')$ is an orbit under
 $(\Homeo 2^\omega)_{\{x_1,\dots,x_n\}}$.  
\end{lem}

\begin{proof}
Homeomorphisms from $(\Homeo 2^\omega)_{\{x_1,\dots,x_n\}}$ preserve types of clo\-pens because they fix $x_1,\dots,x_n$.
Conversely, consider two clopens $b,c$ of the same type $(I,I')$.
If $I=\emptyset$ then taking the union of a homeomorphism $b\to c$ and a homeomorphism $2^\omega\setminus b\to 2^\omega\setminus c$ that fixes $x_1,\dots,x_n$ yields a homeomorphism $\psi \in (\Homeo 2^\omega)_{\{x_1,\dots,x_n\}}$ with $\psi(b)=c$.
An analogous construction deals with the case $I'=\emptyset$.
Consider now the case where $I,I'\neq \emptyset$. Decompose each of $b$ and $c$ according to \eqref{eq:typepart};
let
$b_{ij}$ ($i\in I$, $j\in \N$), $b_{ij}'$ ($i\in I'$, $j\in \N$) be the clopens associated with $b$, and let
$c_{ij}$ ($i\in I$, $j\in \N$), $c_{ij}'$ ($i\in I'$, $j\in \N$) be those for~$c$.
Take arbitrary homeomorphisms $\psi_{ij}\colon b_{ij}\to c_{ij}$ ($i\in I$, $j\in\N$) and
$\psi_{ij}'\colon b_{ij}'\to c_{ij}'$ ($i\in I'$, $j\in\N$).
Their union $\psi$ extended by $\psi(x_i)=x_i$ ($i\in [n]$) is a homeomorphism of $2^\omega$ because of the limit
conditions in \eqref{eq:typepart}, and clearly $\psi(b)=c$, completing the proof.
\end{proof}

\subsection{SIP for the filtered Boolean power}
 We need one more auxiliary result for the proof of Theorem~\ref{thm:SIP}\eqref{it:SIP}.

\begin{lem} \label{lem:invariant}
 Let $\A$ be a finite simple Mal'cev algebra with an idempotent $e_1$,
 let $\B$ be the countable atomless Boolean algebra, let $x_1\in 2^\omega$, let $\D := (\A^\B)^{x_1}_{e_1}$,
   let $h\colon \Aut\D\to(\Homeo 2^\omega)_{x_1}$ be as in Theorem~\ref{thm:AutABxe}\eqref{it:AAB1}, and let
   $g\colon(\Homeo 2^\omega)_{x_1}\to\Aut\D$ be as in~\eqref{eq:defg}.  

   Any subgroup $L$ of $K := \ker h$ with $|K:L| < 2^{\aleph_0}$ that is invariant under conjugation by
   $C := g((\Homeo 2^\omega)_{x_1})$
 contains
\[ K' := \{\varphi\in K \st \varphi_x \in(\Aut\A)_{e_1} \text{ for all } x\in X^\circ \}. \]
\end{lem}   

\begin{proof}
 If $(\Aut\A)_{e_1}$ is trivial, then so is $K'$ and the result clearly holds. So we assume that $(\Aut\A)_{e_1}$
 is non-trivial for the rest of the proof. 
  
 For $\varphi\in K'$, let $\varphi_*\colon X^\circ \to(\Aut\A)_{e_1},\ x \mapsto\varphi_x$, be as in
 Lemma~\ref{la:ppohom}\eqref{it:kerh0}.
 Then $\varphi_*^{-1}(\alpha)$ is clopen in $X^\circ$ for any $\alpha\in(\Aut\A)_{e_1}$ by that
 lemma.

 We will prove the result by showing that
 \[ W := \{\varphi\in K' \st x_1\in\cl( \varphi_*^{-1}(\alpha)) \text{ for all } \alpha\in(\Aut\A)_{e_1} \} \]
 is contained in $L$ and generates $K'$. Note that $\varphi_*^{-1}(\alpha)\neq\emptyset$ for any
 $\varphi\in W$ and $\alpha\in(\Aut\A)_{e_1}$, in particular, $\varphi_*^{-1}(\alpha)$ is clopen of type $([1],[1])$. 
 
 First we claim that
\begin{equation} \label{eq:LcapW}
 L \cap W \neq \emptyset.
\end{equation}   
 To see this, partition $X^\circ$ into a sequence of clopens in $2^\omega$ converging to $x_1$.
 Then subdivide each of the clopens into $|(\Aut\A)_{e_1}|$ non-empty clopens to obtain a partition of
 $X^\circ$ into $|(\Aut\A)_{e_1}|$ disjoint sequences of clopens $b_{i\alpha}$ for $i\in\N,\alpha\in(\Aut\A)_{e_1}$ that converge to $x_1$.
 For any subset $T\subseteq \N$ define $\varphi^T\in K'$ such that
\[
\varphi^T_x := \begin{cases} \alpha & \text{if } x\in b_{i\alpha}\text{ and } i\in T,\\
 \id_A & \text{otherwise}. \end{cases}
\]
Let $\mathcal{T}$ be a family of $2^{\aleph_0}$ subsets of $\N$ such that $T_1\triangle T_2$ (the symmetric difference)
 is infinite for any two distinct $T_1,T_2\in\mathcal{T}$.
 To see that such $\mathcal{T}$ exists, one can consider the equivalence relation on the power-set $2^\N$ defined by
 \[
 T_1\sim T_2\ \Leftrightarrow\ |T_1\Delta T_2|\text{ is finite},
 \]
  note that its equivalence classes are all countable, and then let $\mathcal{~T}$ be a set of representatives of these classes (using AC).
Since $|K:L| < 2^{\aleph_0}$, there exist distinct $T_1,T_2\in\mathcal{T}$ with
 $\varphi:=\varphi^{T_1}(\varphi^{T_2})^{-1}\in L$.  
We claim that $\varphi\in L\cap W$.
  For $i\in T_1\setminus T_2$ and $\alpha\in\Aut \A$ we have
 $\varphi_x^{T_1}=\alpha$  and $\varphi_x^{T_2}=\id_A$ for all $x\in b_{i\alpha}$.
 Similarly, for $i\in T_2\setminus T_1$ and $\alpha\in\Aut \A$ we have
 $\varphi_x^{T_1}=\id_A$  and $\varphi_x^{T_2}=\alpha^{-1}$ for all $x\in b_{i\alpha^{-1}}$.
 It follows that
 \[
 \bigcup_{i\in T_1\setminus T_2} b_{i\alpha}\ \cup \bigcup_{i\in T_2\setminus T_1} b_{i\alpha^{-1}}\subseteq \varphi_*^{-1}(\alpha) \quad \text{for all } \alpha\in\Aut\A.
 \]
 As $T_1,T_2\in\mathcal{T}$, the set $T_1\Delta T_2$ is infinite, and so the set on the left hand side of the above inclusion has $x_1$ in its closure, implying $\varphi \in W$, and completing the proof of the claim that $\varphi\in L\cap W$.

 We show that 
\begin{equation} \label{eq:trans}
 C \text{ acts transitively on } W \text{ by conjugation}.
\end{equation}   
 Let $\sigma,\tau\in W$. For each $\alpha\in(\Aut\A)_{e_1}$ the clopens
 $\sigma_*^{-1}(\alpha)$ and $\tau_*^{-1}(\alpha)$
 have type $([1],[1])$.
 By Lemma~\ref{lem:types} there exists $\psi\in(\Homeo 2^\omega)_{x_1}$ such that
 $\psi(\sigma_*^{-1}(\alpha)) = \tau_*^{-1}(\alpha)$ for all  $\alpha\in(\Aut\A)_{e_1}$.
 Recall that $g(\psi)=\psi^\D\in C$ and $(\tau^{\psi^\D})_x = \tau_{\psi(x)}$. 
 It follows that $\tau^{\psi^\D} = \sigma$ and~\eqref{eq:trans} is proved.

 Since $L$ is normal in $\Aut\D$ by assumption, claims~\eqref{eq:LcapW} and~\eqref{eq:trans} imply $W\subseteq L$.
 It remains to show that
\begin{equation} \label{eq:WK'}
 \langle W \rangle = K'.
\end{equation}   
 For $\alpha\in(\Aut\A)_{e_1}$ and $c$ clopen in $X^\circ$, define $\chi^{c,\alpha}\in K$ by
\[ \chi^{c,\alpha}_x := \begin{cases} \alpha & \text{if } x\in c, \\ \id_A & \text{else.} \end{cases} \]
 To show that $\chi^{c,\alpha}$ is generated by $W$ we distinguish two cases, depending on the type of $c$.
\medskip

\textit{Case 1: $c$ is clopen in $2^\omega$.}
 As in the proof of Claim~\eqref{eq:LcapW} we partition $X^\circ \setminus c$ into $|(\Aut\A)_{e_1}|$ sequences of clopens converging to
 $x_1$. Use these to define $\sigma,\tau\in W$ with $c\subseteq\sigma_*^{-1}(\alpha)$,
 $c\subseteq\tau_*^{-1}(\id)$ and $\sigma_x = \tau_x$ for all $x\in X^{\circ}\setminus c$.
 Then $\chi^{c,\alpha} = \sigma\tau^{-1} \in\langle W\rangle$.
 \medskip
 
 \textit{Case 2: $2^\omega\setminus c$ is clopen in $2^\omega$.}
 This is analogous to Case 1.
 \medskip
 
 \textit{Case 3: $c$ is a disjoint union of clopens in $2^\omega$ converging to $x_1$.}
 This time we partition $c$ itself into $|(\Aut\A)_{e_1}|$ disjoint sequences of clopens $b_{i\beta}$ on $2^\omega$ for
 $i\in\N,\beta\in(\Aut\A)_{e_1}$ that converge to $x_1$. 
 Define $\sigma,\tau\in W$ by
\[
 \sigma_x := \begin{cases} \alpha\beta & \text{if } x\in b_{i\beta},\\
   \id_A & \text{else}, \end{cases}
 \quad \tau_x := \begin{cases} \beta & \text{if } x\in b_{i\beta},\\
   \id_A & \text{else}. \end{cases}
\]
 Then $\chi^{c,\alpha} = \sigma\tau^{-1} \in\langle W\rangle$.

 Since every element in $K'$ is a product of elements $\chi^{c,\alpha}$ for $c$ clopen in $X^\circ$ and
 $\alpha\in(\Aut\A)_{e_1}$, it follows that $K'$ is generated by $W$.
 This proves~\eqref{eq:WK'} and the lemma. 
\end{proof}

 Finally we are ready to prove the main result of this section.

\begin{proof}[Proof of Theorem~\ref{thm:SIP}\eqref{it:SIP}]
 Let $\D := \ABxe$ and let $H$ be a subgroup of $\Aut\D$ of index $<2^{\aleph_0}$.
 Let $G:=(\Homeo 2^\omega)_{\{x_1,\dots,x_n\}}$.
 By Corollary~\ref{cor:isoreduced} we may assume that $e_1,\dots,e_n$ are in distinct $\Aut\A$-orbits. 
 So, by Theorem~\ref{thm:AutABxe}\eqref{it:AAB1} and~\eqref{eq:defg}, $\Aut\D$ is a semidirect product of $K := \ker h$ and
 $C := g(G)$.
 
  Then $|C:H\cap C| < 2^{\aleph_0}$. By Theorem~\ref{thm:Truss} we have $m\geq n$ and a partition of $2^\omega$
  into clopens
  $b_1,\dots,b_m$ with $x_i\in b_i$ for $i\in [n]$ such that
  \[ C' := \{ \psi^\D \st \psi\in G,\ \psi(b_i) = b_i \text{ for } i\in [m] \} \leq H\cap C. \]
 For $i\in [n]$, we will now restrict to a copy of $\Aut(\A^\B)^{x_i}_{e_i}$ inside $\Aut\D$. Let
\[ C_i := \{ (\psi \cup \id_{X\setminus b_i})^\D \st \psi \in (\Homeo\, b_i)_{x_i} \} \leq C' \]
 and   
\[ K_i := \{ \varphi\in K \st \varphi_x = \id_A \text{ for } x\in X^\circ \setminus b_i \}. \]
Since $b_i$ is homeomorphic to $2^\omega$, Theorem~\ref{thm:AutABxe} yields that
 $\Aut(\A^\B)^{x_i}_{e_i}$ is isomorphic to the semidirect product of the normal subgroup $K_i$ and the complement $C_i$.
 Since $C_i\leq H$, also $H\cap K_i$ is invariant under conjugation by $C_i$.
 Further $|K_i:H\cap K_i| < 2^{\aleph_0}$.
 Applying Lemma~\ref{lem:invariant} for $K_iC_i$ as isomorphic copy of $\Aut(\A^\B)^{x_i}_{e_i}$ and $L := H\cap K_i$
 yields that $H\cap K_i$ contains
 \[ K'_i := \{\varphi\in K \st \varphi_x \in(\Aut\A)_{e_i} \text{ for } x\in b_i\setminus\{x_i\},   \varphi_x =\id_A \text{ for } x\in X^\circ\setminus b_i \}. \]
 It follows that
\[  K' := K'_1 \dots K'_n \leq H\cap K. \]
 We complete the proof by showing that $K'C' \leq H$ contains the stabilizer of finitely many elements in $D$.
 Specifically,
for $a$ in $T := \prod_{i=1}^n \{e_i\} \times A^{m-n}$, we define
\[ f_a\colon 2^\omega\to A,\ x\mapsto a_i \text{ if } x\in b_i, i\in [m],\]  
and claim that 
\begin{equation} \label{eq:stab}
 (\Aut\D)_{\{f_a \st a\in T\}} \leq K'C' \leq H.
\end{equation}  
 Consider an element in the above stabilizer decomposed as $\varphi\psi^\D$
 for $\varphi\in K$, $\psi\in G$.
 We first show that $\psi$ stabilizes each set $b_1,\dots,b_m$, so $\psi^\D\in C'$.
 Let $a\in T$. By~\eqref{eq:defpsiD}
\[ f_a = \varphi\psi^\D(f_a) = \varphi(f_a\psi^{-1}). \]
Evaluating both sides at $\psi(x)=\varphi^\circ\psi(x)$ and applying  Lemma~\ref{lem:hphi}\eqref{it:phix}
on the right hand side gives
 \[ f_a(\psi(x)) = \varphi_{\psi(x)}(f_a(x)) \]
 for $x\in X^\circ$.
 If $x\in b_i$ and $\psi(x)\in b_j$ for $i,j\in [m]$, this yields
\begin{equation} \label{eq:ajphiai}
 a_j = \varphi_{\psi(x)}(a_i).
\end{equation}
 So, for all $a\in T$, $a_j$ is uniquely determined by $a_i$ and conversely.
 Suppose that $i\neq j$. Since $|A|>1$, it follows that $i,j\in [n]$.
 But then~\eqref{eq:ajphiai} yields that $e_i,e_j$ are in the same orbit under $\Aut\A$, contradicting our assumption.
 Thus $i=j$ and $\psi^\D\in C'$.

 Now~\eqref{eq:ajphiai} simplifies to 
\[ a_i = \varphi_{x}(a_i) \text{ for all } x\in b_i\cap X^\circ, i \in [m]. \]
 If $i\in [n]$, we obtain $\varphi_x(e_i) = e_i$ for $x\in b_i$.
 If $i\in [m]\setminus [n]$, we obtain $\varphi_x(a_i) = a_i$ for all $a_i\in A$, i.e., $\varphi_x = \id_A$ for all $x\in b_i$.
 Thus $\varphi\in K'$ and~\eqref{eq:stab} is proved, completing the proof that $\Aut\D$ has the SIP.
\end{proof}

\section{Cofinality and Bergman property} \label{sec:Bergman}

In this section we prove Theorem \ref{thm:SIP} \eqref{it:Bergman}.
Throughout we let  $x_1,\dots,x_n\in 2^\omega$ be distinct and $X^\circ:=2^\omega\setminus \{x_1,\dots,x_n\}$. For $\sigma\in G$ let
 \[ \supp\sigma := \{ x\in 2^\omega \st \sigma(x)\neq x \}. \]

\begin{lem} \label{lem:boundconj}
 Let $\sigma,\tau\in(\Homeo 2^\omega)_{\{x_1,\dots,x_n\}}$ be such that
$x_1,\dots,x_n$ are in the closure of $\supp\sigma$.
Then $\tau$ can be written as a product of $8n$ conjugates of $\sigma^{\pm 1}$
by elements in $(\Homeo X)_{\{x_1,\dots,x_n\}}$.
\end{lem}

\begin{proof}
This is a small modification of the final paragraph of the proof of \cite[Theorem 3.17]{Tr:IPG1}.
Use \cite[Lemma 3.16]{Tr:IPG1} to write $\tau=\tau_n\dots \tau_1$, where $\tau_i$  fixes pointwise a clopen $b_i$ containing 
$\{ x_1,\dots, x_{i-1},x_{i+1},\dots, x_n\}$.
We show that every $\tau_i$ can be written as a product of $8$ conjugates of $\sigma^{\pm 1}$. For ease of notation, we
take $i=1$.
By \cite[Lemma 3.15]{Tr:IPG1}, there exist a clopen set $b$ of $X$ which does not contain $x_2,\dots,x_n$, 
and $\mu\in(\Homeo 2^\omega)_{\{x_1,\dots,x_n\}}$ such that $\supp(\sigma\mu\sigma^{-1}\mu^{-1})$ has $x_1$ in its closure
 and is contained within $b$.
Let $c:=b\cup (2^\omega\setminus b_1)$. This again is a clopen that does not contain $x_2,\dots,x_n$, and contains
$\supp(\sigma\mu\sigma^{-1}\mu^{-1})$; furthermore $\supp\tau_1\subseteq c$.
Both 
$\tau_1$ and $\sigma\mu\sigma^{-1}\mu^{-1}$ restrict to homeomorphisms of $c$ since they fix the complement of $c$ in
$2^\omega$. Since $\sigma\mu\sigma^{-1}\mu^{-1}$ does not fix any clopen around $x_1$ pointwise, it follows by
\cite[Corollary 3.14]{Tr:IPG1} that $\tau_1$ can be written as a product of $4$ conjugates of $\sigma\mu\sigma^{-1}\mu^{-1}$,
i.e. a product of $8$ conjugates of $\sigma^{\pm 1}$.
\end{proof}

\begin{lem} \label{le:Xucsc}
  $(\Homeo 2^\omega)_{\{x_1,\dots,x_n\}}$ has uncountable strong cofinality.
\end{lem} 

\begin{proof}
 For $n=0$ this is part of~\cite[Theorem 3.3]{DG:UCPG}.
  
 Assume $n\geq 1$. 
 Let us call clopens of $X^\circ$ of type $([n],[n])$ \emph{good}.
 Recall that $c\subseteq X^\circ$ is a good clopen if and only if both $c$ and $X^\circ \setminus c$ can be partitioned into $n$ sequences of clopens of $2^\omega$ converging to $x_1,\dots,x_n$ respectively.
 We will prove the lemma by applying~\cite[Theorem~2.8]{DG:UCPG} with:
\begin{itemize}
\item  $G := (\Homeo 2^\omega)_{\{x_1,\dots,x_n\}}$;
\item  $\Omega := X^\circ := 2^\omega\setminus\{x_1,\dots,x_n\}$;
\item $\mathfrak{K}$ the set of all good clopens;
 \item
 $F := \{b_1,\dots,b_{n+2}\}$, a set of good clopens that partition $X^\circ$.
\end{itemize}  

We will follow the path offered by part (I) of \cite[Theorem 2.8]{DG:UCPG} and verify that
\begin{itemize}
\item
 the `piecewise-patching' properties (1), (2), (5), (6), (7), (8) of \cite[p. 337]{DG:UCPG} hold;
\item
$G=E^3$, where $E:=\bigcup_{b\in F} G_{\{b\}}$ and $G_{\{b\}}$ denotes the setwise stabilizer of $b$ in $G$.
\end{itemize}

 (1) $\emptyset,X^\circ\not\in\K$ and $\K$ is closed under complements: This follows from the definition.

 (2) Every good set $c$ splits into a disjoint of union of two good sets: Take the partition of $c$ into
 $n$ sequences of disjoint clopens converging to $x_1,\dots,x_n$, respectively, and split each of these clopens
 into two non-empty clopens. This yields $2n$ sequences of disjoint clopens with two of them converging to $x_i$
 for $i\in [n]$. The unions of $n$ of these sequences converging to $x_1,\dots,x_n$, respectively, yields a good
 clopen again. So $c$ can be split into two good clopens.

 (5) $G$ is transitive on $\K$: This follows from Lemma~\ref{lem:types}.

 (6) Let $c,d_1,d_2,e_1,e_2$ be good such that $c=d_1 \dot{\cup} d_2 = e_1 \dot{\cup} e_2$ and
 $f_i\in G$ with $f_i(d_i) = e_i$ for $i\in [2]$. Then $f_1|_{d_1} \cup f_2|_{d_2} \cup \id_{2^\omega\setminus c}$ is in $G$:
 These partial homeomorphisms can be glued together since 
 both $d_1$ and $d_2$ have $x_1,\dots,x_n$ in their closures.
 
 (7) There are disjoint good  $c_i$ for $i\in\N$ such that for all $f_i\in G_{\{c_i\}}$,
 we have $\bigcup_{i\in\N} f_i|_{c_i} \cup \id_{2^\omega\setminus \bigcup_{i\in\N} c_i}$ in $G$:
 Partition $X^\circ$ into $n$ sequences of clopens of $2^\omega$ converging to $x_1,\dots,x_n$, respectively.
 Next split each of these sequences into countably many subsequences. 
 Taking unions of these sequences yields our countably many good clopens $c_i$.

 (8) There exists $m\in\N$ such that for every good clopen $c$, and  any $f,g\in G$ with
 $\supp f,\supp g \subseteq c$ and $\supp g$ containing some good clopen,  we have
 \[ f = (g'_1)^{h_1} \dots (g'_m)^{h_m} \]
 for some $h_i\in G$ with $\supp h_i \subseteq c$ and $g'_i\in \{g,g^{-1}\}$ for all $i\in [m]$:
 This follows from Lemma~\ref{lem:boundconj} with $m=8n$ since the union of any good clopen $c$ with $\{x_1,\dots,x_n\}$
 is homeomorphic to $2^\omega$.

 To check that $G=E^3$ we proceed via a series of claims, motivated by
\cite[Lemma 2.1]{Ga:GCSP}.

\begin{claim} \label{cl:GbGcGb}
Let $b,c,d$ be pairwise disjoint good clopens such that $X^\circ=b\cup c\cup d$ and let $\sigma\in G$ such that  
 $d\setminus\sigma^{-1}(b)$ is good.
Then $\sigma\in G_bG_cG_b$ (pointwise stabilizers).
\end{claim}

\begin{proofclaim}
We aim to define three homeomorphisms $\sigma_1,\sigma_3\in G_b$, $\sigma_2\in G_c$ such that
$\sigma=\sigma_3\sigma_2\sigma_1$,

 Since $(c\cup d)\cap \sigma^{-1}(b)$ is clopen in $X^\circ$ and $d$ is good,
 we have $\tau_1\in G$ such that
\begin{equation} \label{eq:Gb1}
 f := \tau_1((c\cup d)\cap \sigma^{-1}(b)) \text{ is contained in } d \text{ and } d\setminus f \text{ is good}.
\end{equation}
By assumption, the clopen $d\setminus \sigma^{-1}(b)$ is good; in particular, $x_1,\dots, x_n$ belong to its closure.
 Hence the set $(c\cup d)\setminus \sigma^{-1}(b)$ also has all $x_i$ as limit points.
 The complement of this set contains the good clopen $b$, and hence also has $x_1,\dots, x_n$ among its
 limit points.
 We conclude that  $(c\cup d)\setminus \sigma^{-1}(b)$ is good.
 An analogous argument shows that $(c\cup d)\setminus f\in \mathfrak{K}$.
By Lemma \ref{lem:types} there exists a homeomorphism
\begin{equation*}
\tau_2 \colon (c\cup d)\setminus\sigma^{-1}(b)\rightarrow (c\cup d)\setminus f.
\end{equation*}
We can now define the first of our three target homeomorphisms on $X$:
\begin{equation}
\label{eq:Gb6}
\sigma_1(x) := 
\begin{cases}
x & \text{if } x\in b \cup \{x_1,\dots,x_n\}, \\
\tau_1(x) & \text{if } x\in (c\cup d)\cap \sigma^{-1}(b),\\
\tau_2(x) & \text{if } x\in (c\cup d)\cap \sigma^{-1}(c\cup d).
\end{cases}
\end{equation}
To see that $\sigma\in G$, observe that it is defined piecewise by homeomorphisms on three clopens that partition $X^\circ$,
and that the images of those three clopens, namely $b$, $f$ and $(c\cup d)\setminus f$,
also partition $X^\circ$. That $\sigma\in G_b$ is immediate from \eqref{eq:Gb6}.

Next the clopen
\[
(b\cap \sigma^{-1}(c\cup d))\cup (d\setminus f)
\]
is good: indeed, $x_1,\dots,x_n$ are in its closure
because $d\setminus f$ is good by \eqref{eq:Gb1}, while its complement also has $x_1,\dots,x_n$ in its closure
because it contains $c\in\mathfrak{K}$.
Pick a homeomorphism
\begin{equation*}
\tau_3\colon \bigl(b\cap \sigma^{-1}(c\cup d)\bigr)\cup (d\setminus f)\rightarrow d.
\end{equation*}
 Define our second target homeomorphism
 by
 \begin{equation}
 \label{eq:Gb11}
 \sigma_2(x):= 
 \begin{cases}
 x & \text{if } x\in c \cup \{x_1,\dots,x_n\},\\
 \sigma\tau_1^{-1}(x) & \text{if } x\in f,\\
 \sigma(x) & \text{if } x\in b\cap \sigma^{-1}(b),\\
 \tau_3(x) & \text{if } x\in \bigl(b\cap \sigma^{-1}(c\cup d)\bigr)\cup (d\setminus f).
 \end{cases}
 \end{equation}
 Then $\sigma_2\in G_c$ since it is defined piecewise by homeomorphisms on clopens partitioning $X^\circ$ with
 images $c, \sigma(c\cup d)\cap b, \sigma(b)\cap b, d$ also partitioning $X^\circ$.
 
 Finally, define
 \begin{equation}
 \label{eq:Gb13}
 \sigma_3(x) := 
 \begin{cases}
 x & \text{if } x\in b \cup \{x_1,\dots,x_n\},\\
 \sigma\tau_2^{-1}(x) & \text{if } x\in c,\\
 \sigma\tau_3^{-1}(x) & \text{if } x\in \tau_3(b\cap \sigma^{-1}(c\cup d)),\\
 \sigma\tau_2^{-1}\tau_3^{-1}(x) & \text{if } x\in \tau_3(d\setminus f).
 \end{cases}
 \end{equation}
 That $\sigma_3\in G_b$ follows as usual since the images of its constituting partial homeomorphisms
 $b, \sigma\tau_2^{-1}(c), \sigma(b)\cap (c\cup d), \sigma\tau_2^{-1}(d\setminus f)$ partition $X^\circ$.
 
 Now we verify that $\sigma=\sigma_3\sigma_2\sigma_1$ using the definitions \eqref{eq:Gb6}, \eqref{eq:Gb11}, \eqref{eq:Gb13}
 of $\sigma_1,\sigma_2,\sigma_3$ respectively.
 For $x\in X^\circ$
 \begin{align*}
 &\sigma_3\sigma_2\sigma_1(x)
 \\
 =&
 \sigma_3\sigma_2\left( 
 \begin{array}{ll}
 x & \text{if } x\in b\\
 \tau_1(x) & \text{if } x\in (c\cup d)\cap \sigma^{-1}(b)\\
 \tau_2(x) & \text{if } x\in (c\cup d)\cap \sigma^{-1}(c\cup d)
 \end{array}
 \right)
 \\
 =&
 \sigma_3\left(\begin{array}{ll}
 \sigma(x) & \text{if } x\in b\cap \sigma^{-1}(b)\\
 \tau_3(x) & \text{if } x\in b\cap \sigma^{-1}(c\cup d)\\
 \sigma(x) & \text{if } x\in (c\cup d)\cap \sigma^{-1}(b)\\
 \tau_2(x) & \text{if } x\in (c\cup d)\cap \sigma^{-1}(c\cup d) \cap \tau_2^{-1}(c)\\
 \tau_3\tau_2(x) &\text{if } x\in (c\cup d)\cap \sigma^{-1}(c\cup d)\cap \tau_2^{-1}(d\setminus f)
 \end{array}\right)
 \\
 =&
 \begin{cases}
 \sigma(x) & \text{if } x\in b\cap \sigma^{-1}(b)\\
 \sigma\tau_3^{-1}\tau_3(x) & \text{if } x\in b\cap \sigma^{-1}(c\cup d)\\
 \sigma(x) & \text{if } x\in (c\cup d)\cap \sigma^{-1}(b)\\
 \sigma\tau_2^{-1}\tau_2(x) & \text{if } x\in (c\cup d)\cap \sigma^{-1}(c\cup d)\cap \tau_2^{-1}(c)\\
 \sigma\tau_2^{-1}\tau_3^{-1}\tau_3\tau_2(x) & \text{if } 
 x\in (c\cup d)\cap \sigma^{-1}(c\cup d)\cap \tau_2^{-1}(d\setminus f)
 \end{cases}
 \\
 =& \sigma(x).
 \end{align*}
Thus Claim~\ref{cl:GbGcGb} is proved.
\end{proofclaim}

\begin{claim} \label{cl:abutone}
  Let $b,c,d$ be clopens in $X^\circ$ and $k\in [n]$ such that $b\cap c=\emptyset$ and $x_k\in\cl(d)$.
  If $x_k\not\in \cl(d\setminus b)$, then $x_k\in\cl(d\setminus c)$.
\end{claim}

\begin{proofclaim}
 The assumption that $x_k\not\in \cl(d\setminus b)$ means that there exists a clopen neighborhood
 $f$ of $x_k$ in $2^\omega$ such that $f\setminus\{ x_k\}$ is contained in the complement of $d\setminus b$.
 Since $f\cap d$ is clopen in $X^\circ$ and is contained in $d$, this implies $f\cap d\subseteq b$.
Since $b\cap c=\emptyset$, it follows that $(f\cap d)\cap c=\emptyset$.
 Hence $f\cap d\subseteq d\setminus c$ and $x_k\in\cl(d\setminus c)$, as required.
\end{proofclaim}

Recall that $F=\{ b_1,\dots, b_{n+2}\}$ is a set of disjoint good clopens that partition $X^\circ$.

\begin{claim}\label{cl:pigeon}
For every $\sigma\in G$ there exist $i,j\in [n+2]$ such that $\sigma\in G_{b_i}G_{b_j}G_{b_i}$.
\end{claim}

\begin{proofclaim}
 Let $k\in [n]$. By Claim \ref{cl:abutone}, $x_k$ is in the closure of $b_{n+2}\setminus \sigma^{-1} (b_i)$
 for all but possibly one $i\in [n+1]$.
 So, by the pigeonhole principle, there exists $i\in [n+1]$ such that $b_{n+2}\setminus \sigma^{-1}(b_i)$ has all of
 $x_1,\dots,x_n$ in its closure.
 Let $b:= b_i$, $c:=\bigcup_{j\in [n+1]\setminus\{i\}} b_j$ and $d:=b_{n+2}$.
 Then $d\setminus\sigma^{-1}(b)$ is good and $\sigma\in G_bG_cG_b$ by Claim \ref{cl:GbGcGb}.
 For any $j\in [n+1]\setminus\{i\}$ we have $G_c \subseteq G_{b_j}$ and hence $\sigma\in G_{b_i} G_{b_j} G_{b_i}$
 as required.
\end{proofclaim}

To complete the proof of Lemma \ref{le:Xucsc}, recall that $E=\bigcup_{b\in F} G_{\{b\}}$.
Using Claim \ref{cl:pigeon}, for any $\sigma\in G$ there exist $b_i,b_j\in F$ such that
$\sigma\in G_{b_i} G_{b_j} G_{b_i}\subseteq G_{\{b_i\}} G_{\{b_j\}} G_{\{b_i\}}\subseteq E^3$.

 Since $G = E^3$ and all assumptions of case (I) of~\cite[Theorem 2.8]{DG:UCPG} are satisfied,
 $G$ has uncountable strong cofinality.
\end{proof}

 Using the semidirect decomposition of $\Aut\ABxe$ from Theorem~\ref{thm:AutABxe} and the previous lemma,
 it is not hard to show our last main result.

\begin{proof}[Proof of Theorem~\ref{thm:SIP}\eqref{it:Bergman}]
 By Corollary~\ref{cor:isoreduced} we may assume that $e_1,\dots,e_n$ are in distinct $\Aut\A$-orbits. 
 So, by Theorem~\ref{thm:AutABxe}\eqref{it:AAB1} and~\eqref{eq:defg},
 $G := \Aut\ABxe$ is a semidirect product of the groups $K := \ker h$ and
 $C := g((\Homeo X)_{\{x_1,\dots,x_n\}})$.
   
 Seeking a contradiction suppose that $G$ has countable strong cofinality.
 That is, there exists a sequence
\[ H_1 \subseteq H_2 \subseteq \dots \]
 of proper subsets of $G$ such that $H_i^{-1} = H_i$, $H_iH_i \subseteq H_{i+1}$ for all $i\in\N$ and
\[ \bigcup_{i\in\N} H_i = G. \]
 Since $C$ has uncountable strong cofinality by Lemma~\ref{le:Xucsc}, we may assume that $C \subseteq H_1$.

 For $\alpha\in\Aut\A$ and $c$ clopen in $X^\circ$ of type $(I,I')$, assume that there exists $\chi^{c,\alpha}\in K$
 such that 
\[ \chi^{c,\alpha}_x := \begin{cases} \alpha & \text{if } x\in c, \\ \id_A & \text{else.} \end{cases} \]
 Then $\alpha(e_i) = e_i$ for all $i\in I$ by Theorem \ref{thm:AutABxe}\eqref{it:AAB2} and $I\cup I' = [n]$.
 We call such a triple $(\alpha,I, I')$ \emph{legal} and call functions of the form $\chi^{c,\alpha}$ \emph{characteristic}.
 Note that conversely for every legal triple $(\alpha,I, I')$ there exists a clopen $c$ in $X^\circ$ of type $(I,I')$
 and a characteristic $\chi^{c,\alpha}$ in $K$.

 We claim that  any two characteristic functions $\chi^{c_1,\alpha}$, $\chi^{c_2,\alpha}$
 corresponding to the same legal triple $(\alpha,I, I')$ are conjugate.
 Indeed, Lemma~\ref{lem:types} implies that there exists
$\psi\in(\Homeo 2^\omega)_{x_1,\dots,x_n}$ such that $\psi(c_1)=\psi(c_2)$.
A straightforward calculation shows that the automorphism $\psi^\D=g(\psi)\in  C$ associated with $\psi$ conjugates $\chi^{c_2,\alpha}$ to $\chi^{c_1,\alpha}$, proving the claim.
 
 Since there are only finitely many legal triples, 
 we may assume that $H_1$ contains at least one characteristic function $\chi^{c,\alpha}$ for any legal triple.
 Since any two characteristic functions corresponding to the same legal triple are conjugate under $C$ it follows that
  $H_3$ contains all characteristic functions from $K$.
 Since every element of $K$ is a product of $|\Aut\A|$ many characteristic functions from $K$,
 it follows that $K \subseteq H_{|\Aut\A|+3}$.
 Then $G = H_{|\Aut\A|+4}$ contradicting our assumption. Thus $G$ has uncountable strong cofinality.
\end{proof}

\section*{Acknowledgment}

We thank the anonymous referees for their diligent reading and constructive comments which improved the presentation of
this paper.


\begin{thebibliography}{99}
\bibitem{An:ASC}
R.D. Anderson.
\newblock The algebraic simplicity of certain groups of homeomorphisms.
\newblock {\em Amer. J. Math.}, 80:955--963, 1958.

\bibitem{Ap:BPG}
A.B. Apps.
\newblock Boolean powers of groups.
\newblock {\em Math. Proc. Cambridge Philos. Soc.}, 91(3):375--395, 1982.


\bibitem{Ap:ACG}
A.B. Apps, On the structure of $\aleph _{0}$-categorical groups,
{\it J. Algebra} 81:320--339, 1983.



\bibitem{AK:TRA}
R.F. Arens and I.~Kaplansky.
\newblock Topological representation of algebras.
\newblock {\em Trans. Amer. Math. Soc.}, 63:457--481, 1948.

\bibitem{BE:SIP}
R.M. Bryant and D.M. Evans.
\newblock The small index property for free groups and relatively free groups.
\newblock {\em J. London Math. Soc. (2)}, 55(2):363--369, 1997.

\bibitem{BG:ARFG}
R.M. Bryant and J.R.J. Groves.
\newblock On automorphisms of relatively free groups.
\newblock {\em J. Algebra}, 137(1):195--205, 1991.

\bibitem{Bu:BP}
S.~Burris.
\newblock Boolean powers.
\newblock {\em Algebra Universalis}, 5(3):341--360, 1975.

\bibitem{BS:CUA}
S.~Burris and H.P. Sankappanavar.
\newblock {\em A course in universal algebra}.
\newblock Springer, New York Heidelberg Berlin, 1981. \newblock Available from

  \verb+www.math.uwaterloo.ca/~snburris/htdocs/UALG/univ-algebra2012.pdf+.

\bibitem{BW:SCT}
S.~Burris and H. Werner.
\newblock Sheaf constructions and their elementary properties.
\newblock {\em Trans. Amer. Math. Soc.}, 248(2):269--309, 1979.


  
\bibitem{DG:UCPG}
 M. Droste and R. G\"{o}bel.
 Uncountable cofinalities of permutation groups.
{\em J. London Math. Soc.}, 71:335--344, 2005.

\bibitem{DH:GAG}
M. Droste and C. Holland.
\newblock Generating automorphism groups of chains.
\newblock {\em Forum Math.}, 17(4):699--710, 2005.

\bibitem{Ev:EAC}
D.M. Evans.
\newblock Examples of {$\aleph_0$}-categorical structures.
\newblock In {\em Automorphisms of first-order structures}, Oxford Sci. Publ.,
  pages 33--72. Oxford Univ. Press, New York, 1994.

\bibitem{Fo:GBT1}
 A.L. Foster.
 \newblock Generalized ``{B}oolean'' theory of universal algebras. {I}.
   {S}ubdirect sums and normal representation theorem.
 \newblock {\em Math. Z.}, 58:306--336, 1953.

\bibitem{Fo:GBT2}
A.L. Foster.
\newblock Generalized ``{B}oolean'' theory of universal algebras. {II}.
  {I}dentities and subdirect sums of functionally complete algebras.
\newblock {\em Math. Z.}, 59:191--199, 1953.

\bibitem{FM:CTC}
R.~Freese and R.N. McKenzie.
\newblock {\em Commutator theory for congruence modular varieties}, volume 125
  of {\em London Math. Soc. Lecture Note Ser.}
\newblock Cambridge University Press, 1987.
\newblock Available from
  \verb+math.hawaii.edu/~ralph/Commutator/comm.pdf+.

\bibitem{FMMT:ALV3}
R.S. Freese, R.N. McKenzie, G.F. McNulty, and W.F. Taylor.
\newblock {\em Algebras, lattices, varieties. {V}ol. {III}}, volume 269 of {\em
  Mathematical Surveys and Monographs}.
\newblock American Mathematical Society, Providence, RI, 2022.
  
  
  \bibitem{Ga:GCSP}
  F. Galvin, Generating countable sets of permutations. 
  {\em J. London Math. Soc.}, 51:230--242, 1995.
  
  \bibitem{Ha:TG}
  M. Hall Jr., {\it The Theory of Groups}, AMS, Chelsea Publishing, 1999.

\bibitem{HiPl:BPDP}
K. Hickin and J.M. Plotkin, Boolean powers: direct decomposition and isomorphism
types. {\em Trans. Amer. Math. Soc.}, 265:607--621, 1981.

\bibitem{Ho:MT}
W.~Hodges.
\newblock {\em Model theory}, volume~42 of {\em Encyclopedia of Mathematics and
  its Applications}.
\newblock Cambridge University Press, Cambridge, 1993.

\bibitem{HHLS:SIP}
W.~Hodges, I.~Hodkinson, D.~Lascar, and S.~Shelah.
\newblock The small index property for {$\omega$}-stable {$\omega$}-categorical
  structures and for the random graph.
\newblock {\em J. London Math. Soc. (2)}, 48(2):204--218, 1993.

\bibitem{KR:TAG}
A.S. Kechris and C.~Rosendal.
\newblock Turbulence, amalgamation, and generic automorphisms of homogeneous
  structures.
\newblock {\em Proc. Lond. Math. Soc. (3)}, 94(2):302--350, 2007.

\bibitem{Mo:HBA1}
S. Koppelberg,
\newblock {\em Handbook of {B}oolean algebras. {V}ol. 1}.
\newblock North-Holland Publishing Co., Amsterdam, 1989.

\bibitem{Kw:GHC}
A.~Kwiatkowska.
\newblock The group of homeomorphisms of the {C}antor set has ample generics.
\newblock {\em Bull. Lond. Math. Soc.}, 44(6):1132--1146, 2012.

\bibitem{MR:CRN}
A.~Macintyre and J.G. Rosenstein.
\newblock {$\aleph \sb{0}$}-categoricity for rings without nilpotent elements and for {B}oolean structures.
\newblock {\em J. Algebra}, 43(1):129--154, 1976.

\bibitem{Mac:SHS}
D.~Macpherson.
\newblock A survey of homogeneous structures.
\newblock {\em Discrete Math.}, 311(15):1599--1634, 2011.

\bibitem{MR:AGBP}
  P.~Mayr and N. Ru{\v s}kuc.
\newblock Ample generics in automorphism groups of Boolean powers of simple algebras.
\newblock {\em arXiv}, 2025. \verb+https://arxiv.org/abs/2509.20121+



\bibitem{MMT:ALV1}
R.N. McKenzie, G.F. McNulty, and W.F. Taylor.
\newblock {\em Algebras, lattices, varieties, Volume {I}}.
\newblock Wadsworth \& Brooks/Cole Advanced Books \& Software, Monterey,
  California, 1987.
  
\bibitem{Ne:VG}
H.~Neumann.
\newblock {\em Varieties of groups}.
\newblock Springer-Verlag New York, Inc., New York, 1967.

\bibitem{Sc:OC}
J.H. Schmerl,
On $\aleph_0$-categoricity of filtered Boolean extensions,
{\it Algebra Universalis},    8(2): 159--161, 1978.


\bibitem{Tr:IPG1}
J.K. Truss,
Infinite permutation groups. I. Products of conjugacy classes.
\newblock {\em J. Algebra}, 120(2):454--493, 1989.

\bibitem{Tr:IPG2}
J.K. Truss.
\newblock Infinite permutation groups. {II}. {S}ubgroups of small index.
\newblock {\em J. Algebra}, 120(2):494--515, 1989.

\end{thebibliography}
\end{document}